\documentclass[12pt]{article}
\usepackage{amsmath}
\usepackage{amssymb}
\usepackage{vatola}

\def\q{\quad}
\def\qq{\qquad}
\def\mod{\pmod}
\def\t{\text}
\def\f{\frac}
\def\e{\equiv}

\def\phq#1{\varphi(q^{#1})}
\def\psq#1{\psi(q^{#1})}
\def\qtq#1{\q\t{#1}\q}
\def\sls#1#2{(\f{#1}{#2})}
 \def\ls#1#2{\big(\f{#1}{#2}\big)}
\def\Ls#1#2{\Big(\f{#1}{#2}\Big)}
\def\sumn{\sum_{n=0}^{\infty}}
\let \pro=\proclaim
\let \endpro=\endproclaim

\begin{document}
 \par\q\newline
\centerline {\bf Ramanujan's theta functions and linear combinations
of four triangular numbers}

$$\q$$
\centerline{Zhi-Hong Sun}
\par\q\newline
\centerline{School of Mathematics and Statistics}
 \centerline{Huaiyin
Normal University} \centerline{Huaian, Jiangsu 223300, P.R. China}
\centerline{Email: zhsun@hytc.edu.cn} \centerline{Homepage:
http://maths.hytc.edu.cn/szh1.htm}

 \abstract{Let $\Bbb Z$ and $\Bbb Z^+$ be the set of integers
 and the set of positive integers, respectively. For
 $a,b,c,d,n\in\Bbb Z^+$
 let $t(a,b,c,d;n)$ be the number of
 representations of $n$ by
 $ax(x+1)/2+by(y+1)/2+cz(z+1)/2+dw(w+1)/2
 $ $(x,y,z,w\in\Bbb Z)$. In this paper, by using
 Ramanujan's theta functions $\varphi(q)$ and $\psi(q)$
 we present many formulas and conjectures on $t(a,b,c,d;n)$.
 \par\q
 }
 \endabstract

\section*{1. Introduction}

\par  Let $\Bbb Z$, $\Bbb Z^+$ and $\Bbb N$ be the set of
integers, the set of positive integers
 and the set of nonnegative integers, respectively,
 and let
 $\Bbb Z^k=\underbrace{\Bbb Z\times
  \Bbb Z\times\cdots\times\Bbb Z}_{k\ times}.$
  For $a_1,a_2,\ldots,
 a_k\in\Bbb Z^+$ $(k\ge 2)$ and $n\in\Bbb N$  set
$$\align &N(a_1,a_2,\ldots,a_k;n)=\big|\{(x_1,\ldots,x_k)\in \Bbb Z^k\ |
\ n=a_1x_1^2+a_2x_2^2+\cdots+a_kx_k^2 \}\big|,
\\&t(a_1,a_2,\ldots,a_k;n)\\&=\Big|\Big\{(x_1,\ldots,x_k)\in \Bbb
Z^k\ \big|\ n\ =a_1\f{x_1(x_1-1)}2+
a_2\f{x_2(x_2-1)}2+\cdots+a_k\f{x_k(x_k-1)}2\Big\}\Big| \endalign$$
 and
$$C(a_1,\ldots,a_k)=\f{i_1(i_1-1)(i_1-2)(i_1-3)}{4!}
+\f{i_1(i_1-1)i_2}2+i_1i_3,$$ where $i_j$ denotes the number of
elements in $\{a_1,\ldots,a_k\}$ which are equal to $j$. For
convenience we also define
$$t(a_1,a_2,\ldots,a_k;n)=N(a_1,a_2,\ldots,a_k;n)=0\qtq{for}n\not\in
\Bbb N.$$
 In 2005 Adiga, Cooper and
Han [ACH] showed that
$$\aligned &t(a_1,a_2,\ldots,a_k;n)\\&=\f{2}{2+C(a_1,\ldots,a_k)}
N(a_1,\ldots,a_k;8n+a_1+\cdots+a_k) \q\t{for $a_1+\cdots+a_k\le
7$}.\endaligned
 \tag 1.1$$ In
2008 Baruah, Cooper and Hirschhorn [BCH] proved that
$$\aligned &t(a_1,a_2,\ldots,a_k;n)
\\&=\f{2}{2+C(a_1,\ldots,a_k)}(N(a_1,\ldots,a_k;8n+8)-N(a_1,\ldots,a_k;2n+2))
\\&\q\q \t{for $a_1+\cdots+a_k=8$}.\endaligned \tag 1.2$$

\par Ramanujan's theta functions $\varphi(q)$ and $\psi(q)$ are defined
by
$$\varphi(q)=\sum_{n=-\infty}^{\infty}q^{n^2}=1+2\sum_{n=1}^{\infty}
q^{n^2}\qtq{and} \psi(q)=\sum_{n=0}^{\infty}q^{n(n+1)/2}\ (|q|<1).$$
It is evident that for positive integers $a_1,\ldots,a_k$ and
$|q|<1$,
$$\align&\sum_{n=0}^{\infty}N(a_1,\ldots,a_k;n)q^{n}=\varphi(q^{a_1})
\cdots\varphi(q^{a_k}),\tag 1.3
\\&\sum_{n=0}^{\infty}t(a_1,\ldots,a_k;n)q^{n}=2^k\psi(q^{a_1})
\cdots \psi(q^{a_k}).\tag 1.4\endalign$$ There are many identities
involving $\varphi(q)$ and $\psi(q)$.  From [BCH, Lemma 4.1] or [Be]
we know that for $|q|<1$,
 $$\align &\psi(q)^2=\varphi(q)\psi(q^2),\tag
 1.5\\&\varphi(q)=\varphi(q^4)+2q\psi(q^8)=\varphi(q^{16})+2q^{4}\psi(q^{32})
  +2q\psi(q^{8}),\tag 1.6
 \\&\varphi(q)^2=\phq 2^2+4q\psq 4^2=\phq 4^2+4q^2\psq 8^2+4q\psq 4^2,\tag 1.7
  \\&\psi(q)\psi(q^3)=\varphi(q^6)\psi(q^4)+q
 \varphi(q^2)\psi(q^{12}).\tag 1.8
 \endalign$$
By [S1, Lemma 2.4],
 $$ \varphi(q)^2=\phq 8^2+4q^4\psq{16}^2+4q^2\psq 8^2+4q\phq{16}\psq 8
 +8q^5\psq 8\psq{32}.
\tag 1.9$$ By [S1, Lemma 2.3], for $|q|<1$ we have
$$\aligned\varphi(q)\phq 3&=
\phq {16}\phq{48}+4q^{16}\psq{32}\psq{96}+2q \phq {48}\psq
8+2q^3\phq{16}\psq{24} \\&\q+6q^4\psq 8\psq{24}+4q^{13}\psq
8\psq{96}+4q^7\psq{24}\psq{32}.\endaligned \tag 1.10$$

\par Let $a,b,c,d,n\in\Bbb Z^+$. From 1859 to 1866 Liouville made
about 90 conjectures on $N(a,b,c,d;n)$ in a series of papers. Most
conjectures of Liouville have been proved. See  Cooper's survey
paper [C], Dickson's historical comments [D] and Williams' book [W].
 Recently, some connections between
$t(a,b,c,d;n)$ and $N(a,b,c,d;8n+a+b+c+d)$ have been found. See
[ACH,BCH,S1,S3,WS2]. More recently Yao [Y] confirmed some
conjectures posed by the author in [S1]. We also note that the
evaluations of $t(a,b,c,d;n)$ ($a+b+c+d\ge 8$) have been given for
some special values of $(a,b,c,d)$. In [C] Cooper determined
$t(a,b,c,d;n)$ for $(a,b,c,d)=(1,3,3,3),$ $(1,2,2,3),\ (1,3,6,6),\
(1,3,4,4),\ (1,1,2,6)$, $(1,3,12,12)$, in [WS1] Wang and Sun
determined $t(a,b,c,d;n)$ for $(a,b,c,d)=(1,2,2,4),$ $(1,2,4,4),\
(1,1,4,4),\ (1,4,4,4)$, $(1,3,3,9)$,
 $(1,1,9,9),\ (1,9,9,$ $9)$, $(1,1,1,9)$, $(1,3,9,9)$, $(1,1,3,9)$,
 in [WS2]  Wang and Sun
determined $t(a,b,c,d;n)$ for $(a,b,c,d)=(1,1,2,8),$ $
(1,1,2,16),(1,2,3,6),\ (1,3,4,12),\ (1,1,$ $3,4),\ (1,1,5,5),\
(1,5,5,5),\ (1,3,3,12),\ (1,1,1,12),$ $(1,1,3,12),(1,3,3,4)$, in
[S1] Sun determined $t(a,b,c,d;n)$ for $(a,b,c,d)=(1,3,3,6),$
$(1,1,8,8),(1,1,4,8)$. In [S1] the author stated that for odd $a$,
$$\align &t(a,a,2a,4b;4n+3a)=4t(a,2a,4a,b;n),
\q t(a,a,6a,4b;4n+3a)=2t(a,a,6a,b;n),
\\&t(a,a,8a,2b;2n)=t(a,2a,2a,b;n),\q
t(a,a,8a,2b;2n+a)=2t(a,4a,4a,b;n).\endalign$$
\par \par In Section 2, using theta function identities we establish
 new general results for $t(a,b,c,d;n)$. Let $a,b,n\in\Bbb Z^+$
  and $k\in\Bbb N$.
 We show that for odd integers $a$ and $b$,
 $$\align &t(a,2a,2a,2b;n)=\f 12N(a,a,4a,2b;8n+5a+2b),
 \\&t(a,3a,3a,2b;n)=N(3a,3a,4a,2b;8n+7a+2b),
 \\&t(a,a,2a,b;n)=2N(a,4a,8a,b;8n+4a+b)\qtq{for}a\e -b\mod
4,
\\&t(2a,2a,3a,b;n)=\f 13N(a,3a,16a,4b;32n+28a+4b)\qtq{for}a\e -b\mod
4,
\\& t(a,6a,6a,b;n)=\f 13N(a,3a,48a,4b;32n+52a+4b)\q\t{for $a\e b\mod
4$},
\\&t(a,3a,b,3b;n)=4N(a,3a,b,3b;2n+a+b)\qtq{for}n\e (a-b)/2\mod 2.\endalign$$
 For $a,b\in\Bbb Z^+$ let $(a,b)$ be the greatest common divisor of
$a$ and $b$. For an odd prime $p$ and $a\in\Bbb Z$ let $\sls ap$ be
the Legendre symbol. Suppose that $a,b,n\in\Bbb Z^+$, $(a,b)=1$,
$b\not\e 0,-a\mod 4$ and there is an odd prime divisor $p$ of $b$
such that $\sls{a(8n+9a)}p=-1$. We prove that
$$t(a,4a,4a,b;n)=\f
12N(a,4a,4a,b;8n+9a+b).$$ In Section 3, using theta function
identities we establish 31 transformation formulas for
$t(a,b,c,d;n)$. As typical examples, for $a,b,c,n\in\Bbb Z^+$ with
$2\nmid a$ we have
$$\align
&t(a,a,2b,2c;2n+a)=2t(a,4a,b,c;n),
\\&t(a,3a,4b,4c;4n+3a)=2t(3a,4a,b,c;n),
\\&t(a,3a,4b,4c;4n+6a)=2t(a,12a,b,c;n),
\\&t(a,7a,2b,2c;2n+a)=t(a,7a,b,c;n),
\\&t(3a,5a,2b,2c;2n+3a)=t(a,15a,b,c;n),
\\&t(a,15a,2b,2c;2n)=t(3a,5a,b,c;n),
\\&t(a,a,6a,8b;8n+6a)=4t(2a,2a,3a,b;n),
\\&t(2a,3a,3a,8b;8n+12a)=4t(a,6a,6a,b;n).
\endalign$$
In Section 4, we completely determine $t(2,3,3,8;n)$,
$t(1,1,6,24;n)$ and $t(1,1,6,8;n)$ for any positive integer $n$.
 In Section 5 we prove some special relations
 between $t(a,b,c,d;n)$ and $N(a,b,c,d;n)$, and pose many
  challenging conjectures based on calculations on Maple.

\section*{2. New general formulas for $t(a,b,c,d;n)$}
\par In this section we present several general formulas for
$t(a,b,c,d;n)$, which were found by calculations on Maple and proved
by using Ramanujan's theta functions.

 \pro{Theorem 2.1} Let $a,b\in\Bbb Z^+$ with $2\nmid
ab$. For $n=0,1,2,\ldots$ we have
$$t(a,2a,2a,2b;n)=\f 12N(a,a,4a,2b;8n+5a+2b).$$
\endpro
Proof. By (1.6) and (1.9),
$$\align&\sum_{n=0}^{\infty}N(a,a,4a,2b;n)q^n
\\&=\phq a^2\phq{4a}\phq{2b}
\\&=(\phq{8a}^2+4q^{4a}\psq{16a}^2+4q^{2a}\psq{8a}^2+4q^a\phq
{16a}\psq{8a}+8q^{5a}\psq{8a}\psq{32a})
\\&\q\times(\phq{16a}+2q^{4a}\psq{32a})(\phq{8b}+2q^{2b}\psq{16b}).
\endalign$$
Collecting the terms of the form $q^{8n+5a+2b}$ we get
$$\align &\sum_{n=0}^{\infty}N(a,a,4a,2b;8n+5a+2b)q^{8n+5a+2b}
\\&=8q^{5a}\psq{8a}\psq{32a}\cdot \phq{16a}\cdot 2q^{2b}\psq{16b}
\\&\q+4q^a\phq{16a}\psq{8a}\cdot 2q^{4a}\psq{32a}\cdot 2q^{2b}\psq{16b}
\\&=32q^{5a+2b}\psq{8a}\psq{16b}\phq{16a}\psq{32a}
=32q^{5a+2b}\psq{8a}\psq{16a}^2\psq{16b}.\endalign$$ Therefore,
$$\align &\sum_{n=0}^{\infty}N(a,a,4a,2b;8n+5a+2b)q^n
\\&=32\psq{a}\psq{2b}\phq{2a}\psq{4a}=2\sum_{n=0}^{\infty}t(a,2a,4a,2b;n)q^n.
\endalign$$
Hence the result follows.
\par\q

\pro{Theorem 2.2} Let $a,b\in\Bbb Z^+$ with $ab\e -1\mod 4$. For
$n\in\Bbb Z^+$ we have
$$t(a,a,2a,b;n)=2N(a,4a,8a,b;8n+4a+b).$$
\endpro
Proof. It is clear that
$$\align &\sumn N(a,4a,8a,b;n)q^n=\phq a\phq b\phq{4a}\phq{8a}
\\&=(\phq{16a}+2q^{4a}\psq{32a}+2q^a\psq{8a})
(\phq{16b}+2q^{4b}\psq{32b}+2q^b\psq{8b})
\\&\times (\phq{16a}+2q^{4a}\psq{32a})\phq{8a}.
\endalign$$
Thus,
$$\align &\sumn N(a,4a,8a,b;8n+4a+b)q^{8n+4a+b}
\\&=2q^b\psq{8b}(\phq{16a}\cdot 2q^{4a}\psq{32a}+2q^{4a}\psq{32a}
\cdot \phq{16a})\phq{8a}\endalign$$ and so
$$\align &\sumn N(a,4a,8a,b;8n+4a+b)q^n
\\&=8\psq{b}\phq{2a}\psq{4a}\phq{a}=8\psq b\psq{2a}^2\phq a
=8\psq a^2\psq{2a}\psq b
\\&=\f 12\sumn t(a,a,2a,b;n)q^n.\endalign$$
This yields the result.
 \par\q\pro{Corollary 2.1} Suppose
$a,b,n\in\Bbb Z^+$ and $ab\e -1\mod 4$. Then
$$3N(a,4a,8a,b;8n+4a+b)
=N(a,a,2a,b;8n+4a+b)-N(a,a,2a,4b;8n+4a+b).$$
\endpro
Proof. From [S3, Theorem 3.1] we know that
$$t(a,a,2a,b;n)=\f 23(N(a,a,2a,b;8n+4a+b)-N(a,a,2a,4b;8n+4a+b)).$$
This together with Theorem 2.2 yields the result.
\par\q
 \pro{Corollary 2.2} Suppose that $n\in\{0,1,2,\ldots\}$. Then
$$N(1,3,4,8;8n+7)=\f 15N(1,1,2,3;8n+7).$$
\endpro
Proof. By Theorem 2.2, $t(1,1,2,3;n)=2N(1,4,8,3;8n+7)$. By [ACH],
$t(1,1,2,3;n)=\f 25N(1,1,2,3;8n+7)$. Thus the result follows.
\par\q \pro{Theorem 2.3} Let $a,b\in\Bbb
Z^+$ with $2\nmid ab$. For $n\in\Bbb Z^+$ we have
$$t(a,3a,3a,2b;n)=N(3a,3a,4a,2b;8n+7a+2b).$$
\endpro
Proof. Clearly,
$$\align &\sumn N(3a,3a,4a,2b;n)q^n=\phq{3a}^2\phq{4a}\phq{2b}
\\&=(\phq{12a}^2+4q^{6a}\psq{24a}^2+4q^{3a}\psq{12a}^2)
(\phq{8b}+2q^{2b}\psq{16b})\phq{4a}.\endalign$$ Thus,
$$\align &\sumn N(3a,3a,4a,2b;4n+3a+2b)q^{4n+3a+2b}
\\&=4q^{3a}\psq{12a}^2\cdot 2q^{2b}\psq{16b}\cdot\phq{4a}
=8q^{3a+2b}\phq{4a}\psq{12a}^2\psq{16b}\endalign$$ and so
$$\align &\sumn N(3a,3a,4a,2b;4n+3a+2b)q^n
=8\phq a\psq{3a}^2\psq{4b}=8\phq a\phq{3a}\psq{6a}\psq{4b}.
\endalign$$
Applying (1.6) and (1.8) we see that
$$\align &\sumn N(3a,3a,4a,2b;4n+3a+2b)q^n
\\&=8(\phq{4a}+2q^a\psq{8a})(\phq{12a}+2q^{3a}\psq{24a})\psq{6a}\psq{4b}
\\&=8(\phq {4a}\phq{12a}+4q^{4a}\psq {8a}\psq{24a}+2q^a\psq
{2a}\psq {6a})\psq{6a}\psq{4b}.\endalign$$ Hence
$$\sumn N(3a,3a,4a,2b;4(2n+a)+3a+2b)q^{2n+a}
=16q^a\psq{2a}\psq{6a}^2\psq{4b}$$ and so
$$\sumn N(3a,3a,4a,2b;8n+7a+2b)q^n=16\psq a\psq{3a}^2\psq{2b}
=\sumn t(a,3a,3a,2b;n)q^n,$$ which gives the result.
\par\q
\pro{Corollary 2.3} Suppose $n\in\Bbb Z^+$. Then
$$N(3,3,4,6;8n+13)=2N(1,3,12,24;8n+13)=\f 25N(1,3,3,6;8n+13).$$
\endpro
Proof. By Theorem 2.3, $t(1,3,3,6;n)=N(3,3,4,6;8n+13)$. By Theorem
2.2, $t(1,3,3,6;n)=2N(1,3,12,24;8n+13)$. By [S1, Theorem 2.10],
$t(1,3,3,6;n)=\f 25N(1,3,3,6;$ $8n+13)$. Thus, the result follows.
\par\q
\pro{Corollary 2.4} Suppose $a,b,n\in\Bbb Z^+$ and $2\nmid ab$. Then
$$\align &N(3a,4a,12a,2b;8n+7a+2b)=\f 12N(3a,3a,4a,2b;8n+7a+2b),\tag i
\\&N(a,3a,3a,8b;8n+7a+2b)\tag ii\\&\q=N(a,3a,3a,2b;8n+7a+2b)
-2N(3a,3a,4a,2b;8n+7a+2b),
\\&N(a,3a,12a,2b;8n+7a+2b)\tag iii\\&\q=N(a,3a,3a,2b;8n+7a+2b)
-\f 32N(3a,3a,4a,2b;8n+7a+2b).\endalign$$
\endpro
Proof. By [S3, Corollary 4.2],
$t(a,3a,3a,2b;n)=2N(4a,12a,3a,2b;8n+7a+2b)$. This together with
Theorem 2.3 proves (i). By [S3,Theorem 3.1] and Theorem 2.3,
$$\align &N(a,3a,3a,2b;8n+7a+2b)-N(a,3a,3a,8b;8n+7a+2b)
\\&=2t(a,3a,3a,2b;n)=2N(3a,3a,4a,2b;8n+7a+2b).\endalign$$
This yields (ii). By [S3, Theorem 4.5] and Theorem 2.3,
$$\align &\f 23(N(a,3a,3a,2b;8n+7a+2b)-N(a,3a,12a,2b;8n+7a+2b))
\\&=t(a,3a,3a,2b;n)=N(3a,3a,4a,2b;8n+7a+2b),\endalign$$
which yields part(iii). Hence the theorem is proved.
\par\q
\pro{Theorem 2.4} Let $a,b,n\in\Bbb Z^+$ with  $2\nmid ab$. For
$n\in\Bbb Z^+$ with $n\e \f{a-b}2\mod 2$ we have
$$t(a,3a,b,3b;n)=4N(a,3a,b,3b;2n+a+b).$$
\endpro
Proof. Clearly,
$$\align &\sum_{n=0}^{\infty}t(a,3a,b,3b;n)q^n
=16\psq a\psq{3a}\psq b\psq{3b}
\\&=16(\phq{6a}\psq{4a}+q^a\phq{2a}\psq{12a})
(\phq{6b}\psq{4b}+q^b\phq{2b}\psq{12b}).\endalign$$ Thus,
$$\align &\sum_{n=0}^{\infty}t(a,3a,b,3b;2n)q^{2n}
\\&=16\phq{6a}\psq{4a}\phq{6b}\psq{4b}+16q^{a+b}\phq{2a}\psq{12a}\phq{2b}\psq{12b}
\endalign$$
and $$ \align &\sum_{n=0}^{\infty}t(a,3a,b,3b;2n+1)q^{2n+1}
\\&=16q^a\phq{2a}\psq{12a}\phq{6b}\psq{4b}
+16q^b\phq{6a}\psq{4a}\phq{2b}\psq{12b}.
\endalign$$
Therefore,
$$\aligned &\sum_{n=0}^{\infty}t(a,3a,b,3b;2n)q^n
\\&=16\phq{3a}\psq{2a}\phq{3b}\psq{2b}+16q^{\f{a+b}2}\phq{a}\psq{6a}\phq{b}\psq{6b}
\endaligned\tag 2.1$$
and $$ \aligned &\sum_{n=0}^{\infty}t(a,3a,b,3b;2n+1)q^n
\\&=16q^{\f{a-1}2}\phq{a}\psq{6a}\phq{3b}\psq{2b}
+16q^{\f{b-1}2}\phq{3a}\psq{2a}\phq{b}\psq{6b}.
\endaligned\tag 2.2$$
On the other hand,
$$\align &\sum_{n=0}^{\infty}N(a,3a,b,3b;n)q^n
=\phq a\phq{3a}\phq b\phq{3b}
\\&=(\phq{4a}+2q^a\psq{8a})(\phq{12a}+2q^{3a}\psq{24a})
\\&\q\times(\phq{4b}+2q^b\psq{8b})(\phq{12b}+2q^{3b}\psq{24b}).
\endalign$$
Collecting the terms of the form $q^{4n+2}$ yields for $a\e -b\mod
4$,
$$\align &\sum_{n=0}^{\infty}N(a,3a,b,3b;4n+2)q^{4n+2}
\\&=2q^a\psq{8a}\cdot \phq{12a}\cdot \phq{4b}\cdot 2q^{3b}\psq{24b}
\\&\q+2q^b\psq{8b}\cdot \phq{12b}\cdot \phq{4a}\cdot 2q^{3a}\psq{24a}
\\&=4q^{a+b+2}(q^{2b-2}\phq{4b}\psq{24b}\psq{8a}\phq{12a}
+q^{2a-2}\phq{4a}\psq{24a}\psq{8b} \phq{12b}).\endalign$$ Hence
$$\aligned &\sum_{n=0}^{\infty}N(a,3a,b,3b;4n+2)q^n
\\&=4q^{\f{a+b}4}(q^{\f{b-1}2}\phq{b}\psq{6b}\psq{2a}\phq{3a}
+q^{\f{a-1}2}\phq{a}\psq{6a}\psq{2b}\phq{3b}).
\endaligned\tag 2.3$$
Similarly, for $a\e b\mod 4$ we have
$$\align&\sum_{n=0}^{\infty}N(a,3a,b,3b;4n+2)q^{4n+2}
\\&=2q^a\psq{8a}\cdot\phq{12a}\cdot 2q^b\psq{8b}\cdot\phq{12b}
+\phq{4a}\cdot 2q^{3a}\psq{24a}\cdot \phq{4b}\cdot 2q^{3b}\psq{24b}
\\&=4q^{a+b}(\phq{12a}\phq{12b}\psq{8a}\psq{8b}+q^{2a+2b}\phq{4a}\phq{4b}
\psq{24a}\psq{24b})\endalign$$ and so
$$\aligned&\sum_{n=0}^{\infty}N(a,3a,b,3b;4n+2)q^n
\\&=4q^{\f{a+b-2}4}(\phq{3a}\phq{3b}\psq{2a}\psq{2b}+q^{\f{a+b}2}\phq{a}\phq{b}
\psq{6a}\psq{6b}).\endaligned\tag 2.4$$ For $a\e -b\mod 4$ combining
(2.2) with (2.3) gives
$$\sum_{n=0}^{\infty}N(a,3a,b,3b;4n+2)q^n=
\f 4{16}q^{\f{a+b}4}\sum_{n=0}^{\infty}t(a,3a,b,3b;2n+1)q^n$$ and so
$$t(a,3a,b,3b;2n+1)=4N(a,3a,b,3b;4n+2+a+b).$$
 For $a\e b\mod 4$ combining
(2.1) with (2.4) gives
$$\sum_{n=0}^{\infty}N(a,3a,b,3b;4n+2)q^n=
\f 4{16}q^{\f{a+b-2}4}\sum_{n=0}^{\infty}t(a,3a,b,3b;2n)q^n$$ and so
$$t(a,3a,b,3b;2n)=4N(a,3a,b,3b;4n+a+b).$$
This completes the proof.
 \par\q \pro{Theorem 2.5} Let $a,b\in\Bbb
Z^+$ with $ab\e -1\mod 4$. For $n\in\Bbb Z^+$ we have
$$t(2a,2a,3a,b;n)=\f 13N(a,3a,16a,4b;32n+28a+4b).$$
\endpro
Proof.  By (1.10),
$$\align &\sum_{n=0}^{\infty}N(a,3a,16a,4b;n)q^n=\phq
a\phq{3a}\phq{16a}\phq {4b}
\\&=(\phq {16a}\phq{48a}+4q^{16a}\psq{32a}\psq{96a}+2q^a \phq {48a}\psq
{8a}+2q^{3a}\phq{16a}\psq{24a} \\&\q+6q^{4a}\psq
{8a}\psq{24a}+4q^{13a}\psq {8a}\psq{96a}+4q^{7a}\psq{24a}\psq{32a})
\\&\q\times \phq{16a}(\phq{16b}+2q^{4b}\psq{32b}).
\endalign$$
Collecting the terms of the form $q^{8n}$ yields
$$\align &\sum_{n=0}^{\infty}N(a,3a,16a,4b;8n)q^{8n}
\\&=\phq{16a}^2\phq{48a}\phq{16b}+4q^{16a}\phq{16a}\psq{32a}\psq{96a}\phq{16b}
\\&\qq+12q^{4a+4b}\phq{16a}\psq{8a}\psq{24a}\psq{32b}
\endalign$$
and so
$$\align &\sum_{n=0}^{\infty}N(a,3a,16a,4b;8n)q^{n}
\\&=\phq{2a}^2\phq{6a}\phq{2b}+4q^{2a}\phq{2a}\psq{4a}\psq{12a}\phq{2b}
\\&\qq+12q^{(a+b)/2}\phq{2a}\psq{a}\psq{3a}\psq{4b}
\\&=\phq{2a}^2\phq{6a}\phq{2b}+4q^{2a}\phq{2a}\psq{4a}\psq{12a}\phq{2b}
\\&\qq+12q^{(a+b)/2}\phq{2a}\psq{4b}(\phq{6a}\psq{4a}+q^a\phq{2a}\psq{12a}).
\endalign$$
Therefore
$$\sum_{n=0}^{\infty}N(a,3a,16a,4b;8(2n+1))q^{2n+1}
=12q^{(a+b)/2+a}\phq{2a}^2\psq{12a}\psq{4b}$$ and so
$$\align&\sum_{n=0}^{\infty}N(a,3a,16a,4b;8(2n+1))q^n
\\&=12q^{(a+b)/4+(a-1)/2}\phq{a}^2\psq{6a}\psq{2b}
\\&=12q^{(a+b)/4+(a-1)/2}(\phq{2a}^2+4q^a\psq{4a}^2)\psq{6a}\psq{2b}.
\endalign$$
This yields
$$\sum_{n=0}^{\infty}N(a,3a,16a,4b;16n+12a+4b)q^n
=12\phq{2a}^2\psq{6a}\psq{2b}+48q^a\psq{4a}^2\psq{6a}\psq{2b}.$$
Hence
$$\align&
\sum_{n=0}^{\infty}N(a,3a,16a,4b;32n+12a+b)q^{2n}=12\phq{2a}^2\psq{6a}\psq{2b},
\\&\sum_{n=0}^{\infty}N(a,3a,16a,4b;16(2n+1)+12a+4b)q^{2n+1}
=48q^a\psq{4a}^2\psq{6a}\psq{2b}\endalign$$ and so
$$\align&
\sum_{n=0}^{\infty}N(a,3a,16a,4b;32n+12a+4b)q^{n}=12\phq{a}^2\psq{3a}\psq{b},
\\&\sum_{n=0}^{\infty}N(a,3a,16a,4b;16(2n+1)+12a+4b)q^n
=48q^{(a-1)/2}\psq{2a}^2\psq{3a}\psq{b}.\endalign$$ Thus,
$$\sum_{n=0}^{\infty}N(a,3a,16a,4b;32n+28a+4b)q^n=48\psq{2a}^2\psq{3a}\psq{b}.$$
This yields
$$N(a,3a,16a,4b;32n+28a+4b)=3t(2a,2a,3a,b;n).$$
Putting all the above together proves the theorem.
\par\q
 \pro{Theorem 2.6} Let $a,b\in\Bbb Z^+$ with
$ab\e 1\mod 4$. For $n\in\Bbb Z^+$ we have
$$t(a,6a,6a,b;n)=\f 13N(a,3a,48a,4b;32n+52a+4b).$$
\endpro
Proof. By (1.6),
$$\align &\sumn N(a,3a,48a,4b;n)q^n=\phq a\phq{3a}\phq{48a}\phq{4b}
\\&=(\phq{4a}+2q^{4a}\psq{8a})(\phq{12a}+2q^{12a}\psq{24a})\phq{48a}\phq{4b}
\\&=(\phq {4a}\phq{12a}+4q^{4a}\psq {8a}\psq{24a}
+2q^a(\phq{12a}\psq {8a}+q^{2a}\phq {4a}\psq{24a}))
\\&\q\times\phq{48a}\phq{4b}.\endalign$$
Thus,
$$\sumn N(a,3a,48a,4b;4n)q^{4n}
=(\phq {4a}\phq{12a}+4q^{4a}\psq {8a}\psq{24a}) \phq{48a}\phq{4b}.$$
Appealing to (1.6) and (1.8) we see that
$$\align &\sumn N(a,3a,48a,4b;4n)q^n
=(\phq {a}\phq{3a}+4q^{a}\psq {2a}\psq{6a}) \phq{12a}\phq{b}
\\&=(\phq {4a}\phq{12a}+4q^{4a}\psq
{8a}\psq{24a}+6q^a(\phq{12a}\psq {8a}+q^{2a}\phq {4a}\psq{24a}))
\\&\q\times\phq{12a}(\phq{4b}+2q^b\psq{8b}).
\endalign$$
Therefore,
$$\sumn N(a,3a,48a,4b;4(4n+a+b))q^{4n+a+b}
=6q^a\phq{12a}\psq{8a}\cdot\phq{12a}\cdot 2q^b\psq{8b}$$ and so
$$\align&\sumn N(a,3a,48a,4b;4(4n+a+b))q^n
\\&=12\phq{3a}^2\psq{2a}\psq{2b}
=12(\phq{6a}^2+4q^{3a}\psq{12a}^2)\psq{2a}\psq{2b}.\endalign$$ It
then follows that
$$\sumn N(a,3a,48a,4b;4(4(2n+3a)+a+b))q^{2n+3a}
=48q^{3a}\psq{12a}^2\psq{2a}\psq{2b}$$ and therefore
$$\sumn N(a,3a,48a,4b;4(4(2n+3a)+a+b))q^n=48\psq a\psq{6a}^2\psq b
=3\sumn t(a,6a,6a,b;n)q^n.$$ This yields the result.
\par\q
\par For $n\in\{0,1,2,\ldots\}$ let
$$r_3(n)=N(1,1,1;n)=|\{(x,y,z)\in\Bbb Z^3\bigm |
n=x^2+y^2+z^2\}|.$$

\pro{Theorem 2.7} For $n\in\Bbb Z^+$ we have
$$\align &t(1,4,4;n)-\f 12N(1,4,4;8n+9)\\&=t(1,4,4;n)-\f 16r_3(8n+9)
\\&=\cases (-1)^{\f{m+1}2}m&\t{if $8n+9=m^2$ for some $m\in\Bbb Z^+$,}
\\0&\t{otherwise.}\endcases\endalign$$
\endpro
Proof. Since
$$\align \sum_{n=0}^{\infty}t(1,1,8;n)q^n&=8\psi(q)^2\psq 8
=8\varphi(q)\psq 2\psq 8
\\&=8(\phq 4+2q\psq 8)\psq 2\psq 8,\endalign$$
we see that
$$\sum_{n=0}^{\infty}t(1,1,8;2n+1)q^{2n+1}=16q\psq 2\psq 8^2$$
and so
$$\sum_{n=0}^{\infty}t(1,1,8;2n+1)q^n=16\psi(q)\psq 4^2
=2\sum_{n=0}^{\infty}t(1,4,4;n)q^n.$$ Thus,
$$t(1,4,4;n)=\f 12t(1,1,8;2n+1).$$
By [S2, Theorem 3.1],
$$t(1,1,8;2n+1)-\f 13r_3(8n+9)=\cases 2(-1)^{\f{m+1}2}m&\t{if
$8n+9=m^2$ for some $m\in\Bbb Z^+$,}
\\0&\t{otherwise.}\endcases$$
Thus,
$$t(1,4,4;n)-\f 16r_3(8n+9)=\cases (-1)^{\f{m+1}2}m&\t{if
$8n+9=m^2$ for some $m\in\Bbb Z^+$,}
\\0&\t{otherwise.}\endcases\tag 2.5$$
On the other hand,
$$\align\sum_{n=0}^{\infty}N(1,4,4;n)q^n&
=\varphi(q)\phq 4^2= (\phq{16}+2q^4\psq{32}+2q\psq 8)(\phq
8^2+4q^4\psq{16}^2).\endalign$$ Thus,
$$\sum_{n=0}^{\infty}N(1,4,4;8n+1)q^{8n+1}=2q\phq 8^2\psq 8$$
and so
$$\sum_{n=0}^{\infty}N(1,4,4;8n+1)q^n=2\varphi(q)^2\psi(q).$$
Hence
$$\sum_{n=0}^{\infty}N(1,4,4;8n+9)q^n=\f{2\varphi(q)^2\psi(q)-2}{q}.\tag 2.6$$
By [S2, (2.13)],
$$\sum_{n=0}^{\infty}r_3(4n+1)q^n=6\varphi(q)\psi(q)^2=6\varphi(q)^2\psq
2=6(\phq 2^2+4q\psq 4^2)\psq 2.$$ Hence
$$\sum_{n=0}^{\infty}r_3(8n+1)q^{2n}=6\phq 2^2\psq 2$$
and so
$$\sum_{n=0}^{\infty}r_3(8n+1)q^n=6\varphi(q)^2\psi(q).$$
It then follows that
$$\sum_{n=0}^{\infty}r_3(8n+9)q^n=\f{6\varphi(q)^2\psi(q)-6}q.\tag
2.7$$ Combining (2.6) with (2.7) gives $r_3(8n+9)=3N(1,4,4;8n+9)$.
Now putting all the above together proves the theorem.
\par\q
\pro{Corollary 2.5} Let $n\in\Bbb Z^+$ with $n\e 1\mod 3$. Then
$$t(1,4,4;n)=\f 12N(1,4,4;8n+9).$$
\endpro

 \pro{Theorem 2.8} Suppose $a,b,n\in\Bbb Z^+$,
$(a,b)=1$ and there is an odd prime divisor $p$ of $b$ such that
$\sls{a(8n+9a)}p=-1$. Then
$$t(a,4a,4a,b;n)=\f
12(N(a,4a,4a,b;8n+9a+b)-N(a,4a,4a,4b;8n+9a+b)).$$
\endpro
Proof. By Theorem 2.7,
$$\align &t(a,4a,4a,b;n)\\&=\sum_{w\in\Bbb Z}t(a,4a,4a;n-bw(w-1)/2)
=\sum_{w\in\Bbb Z}t\Big(1,4,4;\f{n-bw(w-1)/2}a\Big)
\\&=\f 12\sum_{w\in\Bbb Z} N\Big(1,4,4;8\f{n-bw(w-1)/2}a+9\Big)
\\&=\f 12\sum_{w\in\Bbb Z}N(a,4a,4a;8n-4bw(w-1)+9a)
\\&=\f 12\sum_{w\in\Bbb Z}N(a,4a,4a;8n+9a+b-b(2w-1)^2)
\\&=\f 12\sum_{w\in\Bbb Z}(N(a,4a,4a;8n+9a+b-bw^2)-N(a,4a,4a;8n+9a+b-b(2w)^2)
\\&=\f 12(N(a,4a,4a,b;8n+9a+b)-N(a,4a,4a,4b;8n+9a+b)).
\endalign$$
This proves the theorem.

\pro{Theorem 2.9} Suppose $a,b,n\in\Bbb Z^+$, $(a,b)=1$, $b\not\e
0,-a\mod 4$ and there is an odd prime divisor $p$ of $b$ such that
$\sls{a(8n+9a)}p=-1$. Then
$$t(a,4a,4a,b;n)=\f
12N(a,4a,4a,b;8n+9a+b).$$
\endpro

Proof. Suppose $8n+9a+b=ax^2+4ay^2+4az^2+4bw^2$ for some
$x,y,z,w\in\Bbb Z$. Then $ax^2\e 9a+b\e a+b\mod 4$ and so $b\e
a(x^2-1)\e 0,-a\mod 4$.  But, $b\not\e 0,-a\mod 4$. Hence,
$N(a,4a,4a,4b;8n+9a+b)=0$. Now the result follows from Theorem 2.8.

\pro{Corollary 2.6} Suppose $a,b,n\in\Bbb Z^+$, $3\nmid a$, $3\mid
b$, $b\not\e 0,3a\mod {4}$ and $3\mid n-a$. Then
$$t(a,4a,4a,b;n)=\f
12N(a,4a,4a,b;8n+9a+b).$$
\endpro

\pro{Corollary 2.7} Suppose $m,n\in\Bbb Z^+$, $m\e 1,2\mod 4$ and
$n\e 1,3\mod 5$. Then
$$t(1,4,4,5m;n)=\f 12N(1,4,4,5m;8n+5m+9).$$
\endpro

\pro{Theorem 2.10} Suppose $a,b,n\in\Bbb Z^+$, $(a,b)=1$ and there
is an odd prime divisor $p$ of $b$ such that $\sls{a(4n+5a)}p=-1$.
Then
$$t(a,a,8a,b;n)=\f
12(N(a,a,8a,b;8n+10a+b)-N(a,a,8a,4b;8n+10a+b)).$$
\endpro
Proof. By [S2, Theorem 3.1],
$$\align &t(a,a,8a,b;n)\\&=\sum_{w\in\Bbb Z}t(a,a,8a;n-bw(w-1)/2)
=\sum_{w\in\Bbb Z}t\Big(1,1,8;\f{n-bw(w-1)/2}a\Big)
\\&=\f 12\sum_{w\in\Bbb Z} N\Big(1,1,8;8\f{n-bw(w-1)/2}a+10\Big)
\\&=\f 12\sum_{w\in\Bbb Z}N(a,a,8a;8n-4bw(w-1)+10a)
\\&=\f 12\sum_{w\in\Bbb Z}N(a,a,8a;8n+10a+b-b(2w-1)^2)
\\&=\f 12\sum_{w\in\Bbb Z}(N(a,a,8a;8n+10a+b-bw^2)-N(a,a,8a;8n+10a+b-b(2w)^2)
\\&=\f 12(N(a,a,8a,b;8n+10a+b)-N(a,a,8a,4b;8n+10a+b)).
\endalign$$
This proves the theorem.
\par\q
\pro{Theorem 2.11} Suppose $a,b,n\in\Bbb Z^+$, $(a,b)=1$ and there
is an odd prime divisor $p$ of $b$ such that $\sls{a(4n+5a)}p=-1$.
Assume that $a$ is even or $ab\e 1,4,5\mod 8$ for odd $a$. Then
$$t(a,a,8a,b;n)=\f 12N(a,a,8a,b;8n+10a+b).$$
\endpro
Proof. Suppose $8n+10a+b=ax^2+ay^2+8az^2+4bw^2$ for some
$x,y,z,w\in\Bbb Z$. Then $a(x^2+y^2)\e 2a+b\mod 4$. If $2\mid a$,
then $2\nmid b$ and so $a(x^2+y^2)\not\e 2a+b\mod 4$. We get a
contradiction. If $4\mid ab-1$, then $a\e b\mod 4$ and so
$a(x^2+y^2)\e 3a\mod 4$. Hence $x^2+y^2\e 3\mod 4$. This is
impossible since $t^2\e 0,1\mod 4$ for $t\in\Bbb Z$. If $2\nmid a$
and $ab\e 4\mod 8$, then $2\mid b$ and so $a(x^2+y^2)\e 2a+b\mod 8$.
This yields $x^2+y^2\e 2a^2+ab\e ab+2\e 6\mod 8$. Since $t^2\e
0,1,4\mod 8$ for $t\in\Bbb Z$ we see that $x^2+y^2\not\e 6\mod 8$
and get a contradiction. By the above, $N(a,a,8a,4b;8n+10a+b)=0$.
Now the result follows from Theorem 2.10.
\par\q
\par{\bf Remark 2.1} Theorem 2.11 is a generalization of [S2, Theorem
3.2].

\section*{3. New transformation formulas for $t(a,b,c,d;n)$}
\par\q In this section we present 31 transformation formulas for
$t(a,b,c,d;n)$. \pro{Theorem 3.1} Let $a,b,c\in\Bbb Z^+$ and $2\nmid
a$. For $n\in\Bbb Z^+$ we have
$$\align
&t(a,a,2b,2c;2n+a)=2t(a,4a,b,c;n),\tag 3.1
\\&t(a,3a,4a,2b;2n+a)=t(a,a,6a,b;n),\tag 3.2
\\&t(a,3a,12a,2b;2n)=t(2a,3a,3a,b;n).\tag 3.3
\endalign$$
\endpro
Proof.  Since
$$\align &\sumn t(a,a,2b,2c;n)q^n=16\psq a^2\psq{2b}\psq{2c}
\\&=16\phq a\psq{2a}\psq{2b}\psq{2c}=16(\phq{4a}+2q^a\psq{8a})\psq{2a}
\psq{2b}\psq{2c}\endalign$$ we see that $$\sumn
t(a,a,2b,2c;2n+a)q^{2n+a}=32q^a\psq{8a})\psq{2a} \psq{2b}\psq{2c}$$
and so
$$\sumn t(a,a,2b,2c;2n+a)q^n=32\psq a\psq{4a}\psq b\psq c=2\sumn
t(a,4a,b,c;n)q^n,$$ which yields $t(a,a,2b,2c;2n+a)=2t(a,4a,b,c;n)$.
\par It is clearly that
$$\align &\sumn t(a,3a,4a,2b;n)q^n
\\&=16\psq a\psq{3a}\psq{4a}\psq{2b}
=16(\phq{6a}\psq{4a}+q^a\phq{2a}\psq{12a})\psq{4a}\psq{2b}.
\endalign$$ Thus,
$$\sumn t(a,3a,4a,2b;2n+a)q^{2n+a}
=16q^a\phq{2a}\psq{4a}\psq{12a}\psq{2b}=16q^a\psq{2a}^2\psq{12a}\psq{2b}$$
and so
$$\sumn t(a,3a,4a,2b;2n+a)q^n=16\psq a^2\psq{6a}\psq b
=\sumn t(a,a,6a,b;n)q^n.$$ This yields
$t(a,3a,4a,2b;2n+a)=t(a,a,6a,b;n)$.
\par Note that
$$\align &\sumn t(a,3a,12a,2b;n)q^n
\\&=16\psq a\psq{3a}\psq{12a}\psq{2b}
=16(\phq{6a}\psq{4a}+q^a\phq{2a}\psq{12a})\psq{12a}\psq{2b}.
\endalign$$
We have
$$\sumn t(a,3a,12a,2b;2n)q^{2n}=16\phq{6a}\psq{4a}
\psq{12a}\psq{2b}=16\psq{6a}^2\psq{4a}\psq{2b}$$ and so
$$\sumn t(a,3a,12a,2b;2n)q^n=16\psq {2a}\psq{3a}^2\psq b
=\sumn t(2a,3a,3a,b;n)q^n,$$ which gives $t(a,3a,12a,2b;2n)
=t(2a,3a,3a,b;n)$. This completes the proof.
\par\q
\pro{Theorem 3.2} Suppose $a,b,c\in\Bbb Z^+$ and $2\nmid a$. For
$n\in\Bbb Z^+$ we have
$$\align
&t(a,3a,4b,4c;4n+3a)=2t(3a,4a,b,c;n),\tag 3.4
\\&t(a,3a,4b,4c;4n+6a)=2t(a,12a,b,c;n),\tag 3.5
\\&t(a,3a,48a,4b;4n)=t(a,6a,6a,b;n),\tag 3.6
\\&t(a,3a,16a,4b;4n+a)=t(2a,2a,3a,b;n),\tag 3.7
\\&t(2a,3a,3a,4b;4n+3a)=2t(2a,3a,3a,b;n).\tag 3.8
\endalign$$
\endpro
Proof. It is clear that for $|q|<1$,
$$\align &\sumn t(a,3a,4b,4c;n)q^n
\\&=16\psq a\psq{3a}\psq{4b}\psq{4c}
=16(\phq{6a}\psq{4a}+q^a\phq{2a}\psq{12a})\psq{4b}\psq{4c}
\\&=16((\phq{24a}+2q^{6a}\psq{48a})\psq{4a}+q^a(\phq{8a}
+2q^{2a}\psq{16a})\psq{12a})\psq{4b}\psq{4c}.
\endalign$$
Thus,
$$\align &\sum_{n=0}^{\infty}t(a,3a,4b,4c;4n+6a)q^{4n+6a}
=32q^{6a}\psq{48a}\psq{4a}\psq{4b}\psq{4c},
\\&\sum_{n=0}^{\infty}t(a,3a,4b,4c;4n+3a)q^{4n+3a}
=32q^{3a}\psq{16a}\psq{12a}\psq{4b}\psq{4c},
\\&\sum_{n=0}^{\infty}t(a,3a,4b,48a;4n)q^{4n}=
16\phq{24a}\psq{4a}\psq{4b} \psq{48a},
\\&\sumn t(a,3a,4b,16a;4n+a)q^{4n+a}=16q^a\phq{8a}\psq{12a}\psq{4b}\psq{16a}
\endalign$$ and so
$$\align &\sum_{n=0}^{\infty}t(a,3a,4b,4c;4n+6a)q^n
=32\psq{a}\psq{12a}\psq{b}\psq
c=2\sum_{n=0}^{\infty}t(a,12a,b,c;n)q^n,
\\&\sum_{n=0}^{\infty}t(a,3a,4b,4c;4n+3a)q^n
=32\psq{4a}\psq{3a}\psq{b}\psq{c}=2\sum_{n=0}^{\infty}t(3a,4a,b,c;n)q^n,
\\&\sum_{n=0}^{\infty}t(a,3a,4b,48a;4n)q^n=16\phq {6a}\psq{12a}\psq
a\psq b=16\psq a\psq{6a}^2\psq b,
\\&\sumn t(a,3a,4b,16a;4n+a)q^n=16\phq{2a}\psq{3a}\psq{b}\psq{4a}
=16\psq{2a}^2\psq{3a}\psq b.
\endalign$$
This yields (3.4)-(3.7). To prove (3.8), appealing to Theorem 3.1 we
see that
$$t(2a,3a,3a,4b;4n+3a)=2t(a,3a,12a,2b;2n)=2t(2a,3a,3a,b;n).$$

\par\q
\pro{Theorem 3.3} Suppose $a,b\in\Bbb Z^+$ with $2\nmid a$. For
$n\in\Bbb Z^+$ we have
$$\align
\\&t(2a,3a,3a,8b;8n+9a)=4t(3a,3a,4a,b;n),\tag 3.9
\\&t(a,a,6a,8b;8n+13a)=4t(a,a,12a,b;n),\tag 3.10
\\&t(a,a,6a,8b;8n+4a)=2t(a,a,3a,b;n),\tag 3.11
\\&t(a,a,6a,8b;8n+6a)=4t(2a,2a,3a,b;n),\tag 3.12
\\&t(2a,3a,3a,8b;8n+12a)=4t(a,6a,6a,b;n),\tag 3.13
\\&t(2a,3a,3a,8b;8n+6a)=2t(a,3a,3a,b;n).\tag 3.14
\endalign$$
\endpro
Proof. By Theorems 3.1 and 3.2,
$$\align
&t(2a,3a,3a,8b;8n+9a)=2t(a,3a,12a,4b;4n+3a)=4t(3a,3a,4a,b;n),
\\&t(a,a,6a,8b;8n+13a)=2t(a,3a,4a,4b;4n+6a)=4t(a,a,12a,b;n).
\endalign$$
This proves (3.9) and (3.10). Since $\psi(q)^2=\varphi(q)\psq 2$ and
$\varphi(q)=\phq 4+2q\psq 8$ we see that
$$\align&\sum_{n=0}^{\infty}t(a,a,6a,8b;n)q^n
\\&=16\psq a^2\psq{6a}\psq{8b}=16\phq a\psq{2a}\psq{6a}\psq{8b}
\\&=16(\phq{4a}+2q^a\psq{8a})(\phq{12a}\psq{8a}+q^{2a}\phq{4a}\psq{24a})
\psq{8b}
\\&=16(\phq{16a}+2q^{4a}\psq{32a}+2q^a\psq{8a})
\\&\q\times(\phq{48a}\psq{8a}+2q^{12a}\psq{96a}\psq{8a}
+q^{2a}\phq{16a}\psq{24a}+2q^{6a}\psq{24a}\psq{32a})\psq{8b}.
\endalign$$
Thus,
$$\align &\sum_{n=0}^{\infty}t(a,a,6a,8b;8n+6a)q^{8n+6a}
\\&=16\psq{8b}(\phq{16a}\cdot
2q^{6a}\psq{24a}\psq{32a}+2q^{4a}\psq{32a}\cdot
q^{2a}\phq{16a}\psq{24a})
\\&=64q^{6a}\phq{16a}\psq{32a}\psq{24a}\psq{8b}
=64q^{6a}\psq{16a}^2\psq{24a}\psq{8b}\endalign$$ and
$$\align &\sum_{n=0}^{\infty}t(a,a,6a,8b;8n+4a)q^{8n+4a}
\\&=16\cdot
2q^{4a}(\phq{48a}\psq{32a}+q^{8a}\phq{16a}\psq{96a})\psq{8a}\psq{8b}
\\&=32q^{4a}\psq{8a}^2\psq{24a}\psq{8b}.
\endalign$$
It then follows that
$$\align
&\sum_{n=0}^{\infty}t(a,a,6a,8b;8n+6a)q^n =64\psq {2a}^2\psq{3a}\psq
b=4\sum_{n=0}^{\infty}t(2a,2a,3a,b;n)q^n,
\\&\sum_{n=0}^{\infty}t(a,a,6a,8b;8n+4a)q^n =32\psq a^2\psq{3a}\psq
b=2\sum_{n=0}^{\infty}t(a,a,3a,b;n)q^n,
\endalign$$
which yields (3.11) and (3.12). By (1.5), (1.6) and (1.8),
$$\align &\sum_{n=0}^{\infty}t(2a,3a,3a,8b;n)q^n
\\&=16\psq {2a}\psq {3a}^2\psq {8b}=16\phq {3a}\psq {2a}\psq {6a}\psq {8b}
\\&=16(\phq {12a}+2q^{3a}\psq{24a})(\phq{12a}\psq {8a}+q^{2a}
\phq {4a}\psq{24a})\psq {8b}.\endalign$$ Thus
$$\align &\sum_{n=0}^{\infty}t(2a,3a,3a,8b;4n)q^{4n}=16\phq{12a}^2\psq
{8a}\psq{8b},
\\&\sum_{n=0}^{\infty}t(2a,3a,3a,8b;4n+2a)q^{4n+2a}=16q^{2a}\phq {4a}\psq
{8b}\phq{12a}\psq{24a}
\endalign$$
and so
$$\align &\sum_{n=0}^{\infty}t(2a,3a,3a,8b;4n)q^{n}=16\phq{3a}^2\psq
{2a}\psq{2b} \\&\q=16(\phq{6a}^2+4q^{3a}\psq{12a}^2)\psq
{2a}\psq{2b},
\\&\sum_{n=0}^{\infty}t(2a,3a,3a,8b;4n+2a)q^n=16\phq {a}\psq
{2b}\phq{3a}\psq{6a}\\&\q=16(\phq {4a}+2q^a\psq{8a})(\phq
{12a}+2q^{3a}\psq{24a})\psq{6a}\psq{2b}.
\endalign$$
Hence,
$$\align &\sumn t(2a,3a,3a,8b;4(2n+3a))q^{2n+3a}
=64q^{3a}\psq{2a}\psq{12a}^2\psq{2b},
\\&\sumn t(2a,3a,3a,8b;4(2n+a)+2a)q^{2n+a}
\\&\q=32q^a(\phq{12a}\psq{8a}+q^{2a}\phq{4a}\psq{24a})\psq{6a}\psq{2b}
=32q^a\psq{2a}\psq{6a}^2\psq{2b}.
\endalign$$
This yields
$$\align &\sumn t(2a,3a,3a,8b;8n+12a)q^n
=64\psq{a}\psq{6a}^2\psq{b}=4\sumn t(a,6a,6a,b;n)q^n
\\&\sumn t(2a,3a,3a,8b;8n+6a)q^n=32\psq{a}\psq{3a}^2\psq{b}
=2\sumn t(a,3a,3a,b;n)q^n,
\endalign$$
which yields (3.13) and (3.14). The proof is now complete.
\par\q

\pro{Theorem 3.4} Suppose $a,b,c\in\Bbb Z^+$ with $2\nmid a$. For
$n\in\Bbb Z^+$ we have
$$\align &t(a,7a,2b,2c;2n+a)=t(a,7a,b,c;n),\tag 3.15
\\&t(a,7a,8b,8c;8n+10a)=2t(4a,7a,b,c;n),\tag 3.16
\\&t(a,7a,8b,8c;8n+28a)=2t(a,28a,b,c;n),\tag 3.17
\\&t(a,7a,8a,4b;4n+6a)=t(a,a,14a,b;n),\tag 3.18
\\&t(a,7a,56a,4b;4n)=t(2a,7a,7a,b;n).\tag 3.19
\endalign$$
\endpro
Proof. By [Be, p.315],
$$\psi(q)\psq 7=\psq 8\phq{28}+q\psq 2\psq{14}+q^6\phq
4\psq{56}.\tag 3.20$$ Thus,
$$\psq 2\psq {14}=\psq {16}\phq{56}+q^2\psq 4\psq{28}+q^{12}\phq
8\psq{112}.\tag 3.21$$ Combining (3.20) with (3.21) gives
$$\aligned \psi(q)\psq 7&=\psq 8\phq{28}+q^6\phq
4\psq{56}+q\psq{16}\phq{56}\\&\q+q^3\psq 4\psq{28}+q^{13}\phq
8\psq{112}.
\endaligned\tag 3.22$$
Using (3.20) we see that
$$\align&\sumn t(a,7a,2b,2c;n)q^n=16\psq a\psq{7a}\psq{2b}\psq{2c}
\\&\q=16(\psq{8a}\phq{28a}+q^{6a}\phq{4a}\psq{56a}+q^a\psq{2a}\psq{14a})
\psq{2b}\psq{2c}\endalign$$ and so
$$\sumn
t(a,7a,2b,2c;2n+a)q^{2n+a}=16q^a\psq{2a}\psq{14a}\psq{2b}\psq{2c}.$$
This yields
$$\sumn
t(a,7a,2b,2c;2n+a)q^n=16\psq{a}\psq{7a}\psq{b}\psq{c}.$$ Hence
(3.15) is true. Applying (3.22) we see that
$$\aligned&\sumn t(a,7a,4b,4c;n)q^n=16\psq a\psq{7a}\psq{4b}\psq{4c}
\\&=16(\psq {8a}\phq{28a}+q^{6a}\phq
{4a}\psq{56a}+q^a\psq{16a}\phq{56a}\\&\q+q^{3a}\psq
{4a}\psq{28a}+q^{13a}\phq
{8a}\psq{112a})\psq{4b}\psq{4c}.\endaligned$$ Thus,
$$\align &\sumn
t(a,7a,4b,4c;4n)q^{4n}=16\psq{8a}\phq{28a}\psq{4b}\psq{4c},
\\&\sumn
t(a,7a,4b,4c;4n+6a)q^{4n+6a}=16q^{6a}\phq{4a}\psq{56a}\psq{4b}\psq{4c}
\endalign$$
and so
$$\align &\sumn
t(a,7a,4b,4c;4n)q^n=16\psq{2a}\phq{7a}\psq{b}\psq{c},\tag 3.23
\\&\sumn
t(a,7a,4b,4c;4n+6a)q^n=16\phq{a}\psq{14a}\psq{b}\psq{c}.\tag 3.24
\endalign$$
By (3.24) and (1.6),
$$\sumn t(a,7a,8b,8c;4n+6a)q^n
=16(\phq{4a}+2q^a\psq{8a})\psq{14a}\psq{2b}\psq{2c}.$$ Therefore
$$\sumn t(a,7a,8b,8c;4(2n+a)+6a)q^{2n+a}
=32q^a\psq{8a}\psq{14a}\psq{2b}\psq{2c}$$ and so
$$\sumn t(a,7a,8b,8c;4(2n+a)+6a)q^n
=32\psq{4a}\psq{7a}\psq{b}\psq{c}.$$ This yields (3.16). Similarly,
from (3.23) and (1.6),
$$\sumn
t(a,7a,8b,8c;4n)q^n
=16(\phq{28a}+2q^{7a}\psq{56a})\psq{2a}\psq{2b}\psq{2c}.$$
Therefore,
$$\sumn
t(a,7a,8b,8c;4(2n+7a))q^{2n+7a}
=32q^{7a}\psq{56a}\psq{2a}\psq{2b}\psq{2c}$$ and so
$$\sumn t(a,7a,8b,8c;4(2n+7a))q^n=32\psq a\psq{28a}\psq b\psq c,$$
which gives (3.17).
 Recall that
$\varphi(q)\psq 2=\psi(q)^2$. From (3.23) and (3.24) we see that
$$\align &\sumn
t(a,7a,56a,4b;4n)q^n=16\psq{2a}\phq{7a}\psq{14a}\psq
b=16\psq{2a}\psq{7a}^2\psq b, \\&\sumn
t(a,7a,8a,4b;4n+6a)q^n=16\phq{a}\psq{2a}\psq{14a}\psq{b}=16\psq
a^2\psq{14a}\psq b,\endalign$$ which yields (3.18) and (3.19). The
proof is now complete.
\par\q
\par From [Be, p.377] we know that for $|q|<1$,
$$\psq 3\psq 5=\phq {60}\psq 8+q^3\psq 2\psq{30}+q^{14}\phq
4\psq{120}.\tag 3.25$$ From [XZ, (2.10)] we know that
$$\psi(q)\psq{15}=\psq
6\psq{10}+q\phq{20}\psq{24}+q^3\phq{12}\psq{40}.\tag 3.26$$ Thus,
$$\psi(q^2)\psq{30}=\psq
{12}\psq{20}+q^2\phq{40}\psq{48}+q^6\phq{24}\psq{80}.$$ This
together with (3.25) yields
$$\aligned &\psq 3\psq 5=\phq{60}\psq
8+q^{14}\phq 4\psq{120}+q^3\psq{12}\psq{20}
\\&\qq\qq\qq+q^5\phq{40}\psq{48}+q^9\phq{24}\psq{80}.\endaligned\tag
3.27$$  By (3.25), $$ \psq 6\psq {10}=\phq {120}\psq {16}+q^6\psq
4\psq{60}+q^{28}\phq 8\psq{240}.$$ Combining this with (3.26) gives
$$\aligned&\psi(q)\psq{15}=\phq{120}\psq{16}+q^{28}\phq 8\psq{240}+q^6\psq
4\psq{60}\\&\qq\qq\qq+q\phq{20}\psq{24}+q^3\phq{12}\psq{40}.
\endaligned\tag 3.28$$
\pro{Theorem 3.5} Suppose $a,b,c\in\Bbb Z^+$ with $2\nmid a$. For
$n\in\Bbb Z^+$ we have
$$\align &t(3a,5a,2b,2c;2n+3a)=t(a,15a,b,c;n),\tag 3.29
\\&t(a,15a,2b,2c;2n)=t(3a,5a,b,c;n),\tag 3.30
\\&t(3a,5a,4b,4c;4n+3a)=t(3a,5a,b,c;n),\tag 3.31
\\&t(a,15a,4b,4c;4n+6a)=t(a,15a,b,c;n),\tag 3.32
\\&t(3a,5a,8b,8c;8n+18a)=2t(4a,15a,b,c;n),\tag 3.33
\\&t(3a,5a,8b,8c;8n+60a)=2t(a,60a,b,c;n),\tag 3.34
\\&t(3a,5a,8a,4b;4n+14a)=t(a,a,30a,b;n),\tag 3.35
\\&t(3a,5a,120a,4b;4n)=t(2a,15a,15a,b;n),\tag 3.36
\\&t(a,15a,8b,8c;8n+15a)=2t(5a,12a,b,c;n),\tag 3.37
\\&t(a,15a,8b,8c;8n+21a)=2t(3a,20a,b,c;n),\tag 3.38
\\&t(a,15a,24a,4b;4n+3a)=t(3a,3a,10a,b;n),\tag 3.39
\\&t(a,15a,40a,4b;4n+a)=t(5a,5a,6a,b;n).\tag 3.40
\endalign$$
\endpro
Proof. Using (3.25) and (3.26) we see that
$$\align &\sumn
t(3a,5a,2b,2c;n)q^n=16\psq{3a}\psq{5a}\psq{2b}\psq{2c}
\\&\q=16(\phq{60a}\psq{8a}+q^{14a}\phq{4a}\psq{120a}
+q^{3a}\psq{2a}\psq{30a})\psq{2b}\psq{2c},
\\&\sumn
t(a,15a,2b,2c;n)q^n=16\psq{a}\psq{15a}\psq{2b}\psq{2c}
\\&\q=16(\psq{6a}\psq{10a}+q^{a}\phq{20a}\psq{24a}
+q^{3a}\phq{12a}\psq{40a})\psq{2b}\psq{2c}.\endalign$$ Hence
$$\align &\sumn
t(3a,5a,2b,2c;2n+3a)q^{2n+3a}=16q^{3a}\psq{2a}\psq{30a}\psq{2b}\psq{2c},
\\&\sumn
t(a,15a,2b,2c;2n)q^{2n}=16\psq{6a}\psq{10a}\psq{2b}\psq{2c}.\endalign$$
It then follows that
$$\align &\sumn
t(3a,5a,2b,2c;2n+3a)q^n=16\psq{a}\psq{15a}\psq{b}\psq{c},
\\&\sumn
t(a,15a,2b,2c;2n)q^n=16\psq{3a}\psq{5a}\psq{b}\psq{c},\endalign$$
which implies (3.29) and (3.30). By (3.29) and (3.30),
$$\align &t(3a,5a,4b,4c;4n+3a)=t(a,15a,2b,2c;2n)=t(3a,5a,b,c;n),
\\&t(a,15a,4b,4c;4n+6a)=t(3a,5a,2b,2c;2n+3a)=t(a,15a,b,c;n).\endalign$$
Thus, (3.31) and (3.32) hold. Appealing to (3.27),
$$\align &\sumn
t(3a,5a,4b,4c;n)q^n=16\psq{3a}\psq{5a}\psq{4b}\psq{4c}
\\&\q=16(\phq{60a}\psq
{8a}+q^{14a}\phq {4a}\psq{120a}+q^{3a}\psq{12a}\psq{20a}
\\&\qq+q^{5a}\phq{40a}\psq{48a}+q^{9a}\phq{24a}\psq{80a})\psq{4b}\psq{4c}.
\endalign$$
Thus,
$$\align&\sumn t(3a,5a,4b,4c;4n)q^{4n}=16\phq{60a}
\psq{8a}\psq{4b}\psq{4c},
\\&\sumn t(3a,5a,4b,4c;4n+14a)q^{4n+14a}=16q^{14a}\phq{4a}
\psq{120a}\psq{4b}\psq{4c},\endalign$$ which yields
$$\align&\sumn t(3a,5a,4b,4c;4n)q^n=16\phq{15a}
\psq{2a}\psq{b}\psq{c},\tag 3.41
\\&\sumn t(3a,5a,4b,4c;4n+14a)q^n=16\phq{a}
\psq{30a}\psq{b}\psq{c}.\tag 3.42\endalign$$ From (3.42), (3.41) and
(1.6) we see that
$$\align &\sumn t(3a,5a,8b,8c;4n+14a)q^n
=16(\phq{4a}+2q^a\psq{8a})\psq{30a}\psq{2b}\psq{2c},
\\ &\sumn
t(3a,5a,8b,8c;4n)q^n
=16(\phq{60a}+2q^{15a}\psq{120a})\psq{2a}\psq{2b}\psq{2c}.
\endalign$$
Thus,
$$\align &\sumn t(3a,5a,8b,8c;4(2n+a)+14a)q^{2n+a}
=32q^a\psq{8a}\psq{30a}\psq{2b}\psq{2c},
\\ &\sumn
t(3a,5a,8b,8c;4(2n+15a))q^{2n+15a}
=32q^{15a}\psq{120a}\psq{2a}\psq{2b}\psq{2c}
\endalign$$
and so
$$\align &\sumn t(3a,5a,8b,8c;4(2n+a)+14a)q^n
=32\psq{4a}\psq{15a}\psq{b}\psq{c},
\\ &\sumn
t(3a,5a,8b,8c;4(2n+15a))q^n =32\psq{a}\psq{60a}\psq{b}\psq{c},
\endalign$$
which yields (3.33) and (3.34). Note that $\varphi(q)\psq
2=\psi(q)^2$. Taking $c=2a$ in (3.42) and $c=30a$ in (3.41) yields
(3.35) and (3.36). Using (3.28) we see that
$$\align &\sumn t(a,15a,4b,4c;n)q^n=16\psq
a\psq{15a}\psq{4b}\psq{4c}
\\&\q=16(\phq{120a}\psq{16a}+q^{28a}\phq {8a}\psq{240a}+q^{6a}\psq
{4a}\psq{60a}\\&\qq+q^a\phq{20a}\psq{24a}+q^{3a}\phq{12a}\psq{40a})
\psq{4b}\psq{4c}.\endalign$$
 From this it follows that
 $$\align &\sumn
t(a,15a,4b,4c;4n+a)q^{4n+a}=16q^a\phq{20a}\psq{24a}\psq{4b}\psq{4c},
\\&\sumn
t(a,15a,4b,4c;4n+3a)q^{4n+3a}=16q^{3a}\phq{12a}\psq{40a}\psq{4b}\psq{4c}
\endalign$$
and so
$$\align &\sumn
t(a,15a,4b,4c;4n+a)q^n=16\phq{5a}\psq{6a}\psq{b}\psq{c},\tag 3.43
\\&\sumn
t(a,15a,4b,4c;4n+3a)q^n=16\phq{3a}\psq{10a}\psq{b}\psq{c}.\tag 3.44
\endalign$$
Hence, applying (1.6) we get
$$\align &\sumn
t(a,15a,8b,8c;4n+a)q^n=16(\phq{20a}+2q^{5a}\psq{40a})
\psq{6a}\psq{2b}\psq{2c},
\\&\sumn
t(a,15a,8b,8c;4n+3a)q^n=16(\phq{12a}+2q^{3a}\psq{24a})
\psq{10a}\psq{2b}\psq{2c}
\endalign$$
and so
$$\align &\sumn
t(a,15a,8b,8c;4(2n+5a)+a)q^{2n+5a}=32q^{5a}\psq{40a}
\psq{6a}\psq{2b}\psq{2c},
\\&\sumn
t(a,15a,8b,8c;4(2n+3a)+3a)q^{2n+3a}=32q^{3a}\psq{24a}
\psq{10a}\psq{2b}\psq{2c}.
\endalign$$
It then follows that
$$\align &\sumn
t(a,15a,8b,8c;4(2n+5a)+a)q^n=32\psq{20a} \psq{3a}\psq{b}\psq{c},
\\&\sumn
t(a,15a,8b,8c;4(2n+3a)+3a)q^n=32\psq{12a} \psq{5a}\psq{b}\psq{c},
\endalign$$
which yields (3.37) and (3.38). Finally, taking $c=6a$ in (3.44) and
$c=10a$ in (3.43) and then applying (1.5) we deduce (3.39) and
(3.40). The proof is now complete.
\par\q
\section*{4. Evaluation of  $t(2,3,3,8;n)$, $t(1,1,6,24;n)$ and $t(1,1,6,8;n)$}
\par In this section we determine  $t(2,3,3,8;n)$, $t(1,1,6,24;n)$
and $t(1,1,6,8;n)$ for any positive integer $n$.
 \pro{Theorem 4.1} Let $n$ be a positive integer.
\par $(\t{\rm i})$ If $2n+5=3^{\beta}n_1\ (3\nmid n_1)$, then
$$t(1,1,6,24;2n+1)=t(2,3,3,8;2n+3)=4(\sigma(n_1)-(-1)^na(2n+5)),$$
where $\{a(n)\}$ is given by
$$q\prod_{k=0}^{\infty}(1-q^{2k})(1-q^{4k})(1-q^{6k})(1-q^{12k})
=\sum_{n=1}^{\infty}a(n)q^n\q(|q|<1).$$
\par $(\t{\rm ii})$ If $n+1=2^{\alpha}3^{\beta}n_1$ with $2\nmid
n_1$ and $3\nmid n_1$, then
$$t(1,1,6,24;2n-2)+t(2,3,3,8;2n)=2^{\alpha+4}\sigma(n_1).$$
\endpro
Proof. By (3.1), $t(1,1,6,24;2n+1)=2t(1,3,4,12;n)=t(2,3,3,8;2n+3)$.
By [WS2, Theorem 3.4], if $2n+5=3^{\beta}n_1$ with $3\nmid n_1$,
then $t(1,3,4,12;n)=2(\sigma(n_1)-(-1)^na(2n+5)).$ Thus (i) is true.
\par Now we prove (ii).
Since $\psi(q)^2=\varphi(q)\psq 2=(\phq 4+2q\psq 8)\psq 2$ we see
that
$$\sum_{n=0}^{\infty}t(1,1,6,24;n)q^n=16\psi(q)^2\psq 6\psq {24}
=(16\phq 4+32q\psq 8)\psq 2\psq 6\psq{24}.$$ Thus,
$$\sum_{n=0}^{\infty}t(1,1,6,24;2n)q^{2n}
=16\phq 4\psq 2\psq 6\psq{24}$$ and so
$$\sum_{n=0}^{\infty}t(1,1,6,24;2n)q^n
=16\phq 2\psi(q)\psq 3\psq{12}.\tag 4.1$$  Similarly,
$$\sum_{n=0}^{\infty}t(2,3,3,8;n)q^n=16\psq 2\psq 3^2\psq 8
=(16\phq {12}+32q^3\psq {24})\psq 6\psq 2\psq 8.$$ Thus,
$$\sum_{n=0}^{\infty}t(2,3,3,8;2n)q^{2n}=16
\psq 2\psq 6 \psq 8\phq{12}$$ and so
$$\sum_{n=0}^{\infty}t(2,3,3,8;2n)q^n=16
\psi(q)\psq 3 \psq 4\phq{6}.\tag 4.2$$
 From (4.1) and (4.2) we see
that
$$\align &\sum_{n=0}^{\infty}(t(2,3,3,8;2n)+t(1,1,6,24;2n-2))q^n
\\&=16\psi(q)\psq 3(\phq{6}\psq 4+q\phq 2\psq {12})
=16\psi(q)^2\psq 3^2.\endalign$$ Hence
$$t(2,3,3,8;2n)+t(1,1,6,24;2n-2)=t(1,1,3,3;n).\tag 4.3$$
By [WS1, Lemma 4.1], if $n+1=2^{\alpha}3^{\beta}n_1$ with $2\nmid
n_1$ and $3\nmid n_1$, then
$$t(1,1,3,3;n)=2^{\alpha+4}\sigma(n_1).\tag 4.4$$
Thus (ii) holds and the proof is complete.
\par\q
\pro{Theorem 4.2} Let $n$
be a positive integer.
\par $(\t{\rm i})$ If $n+1=2^{\alpha}3^{\beta}n_1$ with $2\nmid n_1$
and $3\nmid n_1$, then
$$t(2,3,3,8;4n+2)=t(1,1,6,24;4n)=2^{\alpha+4}\sigma(n_1).$$
\par $(\t{\rm ii})$ If $2n+1=3^{\beta}n_1\ (3\nmid n_1)$, then
$$\align &t(2,3,3,8;4n)=8(\sigma(n_1)+a(2n+1)),
\\&t(1,1,6,24;4n-2)=8(\sigma(n_1)-a(2n+1)).\endalign$$
\endpro
Proof.  If $n+1=2^{\alpha}3^{\beta}n_1$ with $2\nmid n_1$ and
$3\nmid  n_1$, from Theorem 4.1 we have
$$t(1,1,6,24;4n)+t(2,3,3,8;4n+2)=2^{\alpha+5}\sigma(n_1).$$
By [S1, Theorem 2.14], $t(1,1,6,24;4n)=2^{\alpha+4}\sigma(n_1).$
Thus, $t(2,3,3,8;4n+2)=2^{\alpha+4}\sigma(n_1).$ This proves (i).
\par Now we consider (ii). Suppose $2n+1=3^{\beta}n_1\ (3\nmid n_1)$.
We first assume that $n$ is odd and $n=2m+1$. By (3.13),
$t(2,3,3,8;8m+12)=4t(1,1,6,6;m).$ This together with [S3, Theorem
4.15] yields
 $$t(2,3,3,8;8m+4)=4t(1,1,6,6;m-1)=8(\sigma(n_1)+a(4m+3)). $$
Now, appealing to Theorem 4.1(ii) we get
$$t(1,1,6,24;8m+2)=16\sigma(n_1)-t(2,3,3,8;8m+4)
=8(\sigma(n_1)-a(4m+3)).$$ From now on suppose that $n$ is even and
$n=2m$. From (3.12) and [S3, Theorem 4.15] we see that
$$t(1,1,6,24;8m-2)=4t(2,2,3,3;m-1)=8(\sigma(n_1)-a(4m+1)).$$
Now applying Theorem 4.1(ii) gives
$$t(2,3,3,8;8m)=16\sigma(n_1)-t(1,1,6,24;8m-2)=8(\sigma(n_1)+a(4m+1)).$$
 Putting all the above together proves the theorem.
\par\q
\pro{Theorem 4.3} Let $n\in\Bbb Z^+$.
\par $(\t{\rm i})$ If $n+1=2^{\alpha}3^{\beta}n_1$ with $2\nmid n_1$
and $3\nmid n_1$, then
$$t(1,1,6,8;2n)=2^{\alpha+2}\Big(3^{\beta+1}\Ls 3{n_1}+(-1)^{\alpha+\beta+
\f{n_1-1}2}\Big)\sum_{d\mid n_1}d\Ls 3d .$$
\par $(\t{\rm ii})$ If $2n+3=3^{\beta}n_1$ with $3\nmid n_1$, then
$$t(1,1,6,8;2n+1)=2\Big(3^{\beta+1}\Ls 3{n_1}+(-1)^n\Big)
\sum_{d\mid n_1}d\Ls 3d-4\sum\Sb a,b\in\Bbb Z^+, 2\nmid
a\\8n+12=a^2+3b^2\endSb (-1)^{(a-1)/2}a.$$
\endpro
Proof. From [S1, Theorem 2.7],
$$t(1,1,6,8;2n)=t(1,2,2,3;n)\qtq{and}t(1,1,6,8;2n+1)=2t(1,3,4,4;n).$$
By [C, (5.6) and Theorem 5.4], if $n+1=2^{\alpha}3^{\beta}n_1$ with
$2\nmid n_1$ and $3\nmid n_1$, then
$$\align t(1,2,2,3;n)&=2^{\alpha+2}\Big(3^{\beta+1}+(-1)^{\alpha+\beta}\Ls{-3}{n_1}
\Big)\sum_{d\mid n_1}\f {n_1}d\Ls 3d
\\&=2^{\alpha+2}\Big(3^{\beta+1}+(-1)^{\alpha+\beta}\Ls{-3}{n_1}
\Big)\sum_{d\mid n_1}d\Ls 3{n_1/d}
\\&=2^{\alpha+2}\Big(3^{\beta+1}\Ls
3{n_1}+(-1)^{\alpha+\beta+(n_1-1)/2}\Big)\sum_{d\mid n_1}d\Ls 3d.
\endalign$$ Thus part(i) is true.
\par Let us consider part(ii).
 By [C, (5.6) and Theorem 5.6], if
$2n+3=3^{\beta}n_1$ with $3\nmid n_1$, then
$$t(1,3,4,4;n)=\Big(3^{\beta+1}-(-1)^{\beta}\Ls{-3}{n_1}\Big)
\sum_{d\mid n_1}\f{n_1}d\Ls 3d-2g_1(2n+3),$$ where $\{g_1(n)\}$ is
given by
$$q\psq
6\prod_{k=1}^{\infty}(1-q^{2k})^3=\sum_{n=1}^{\infty}g_1(n)q^n.$$
Using Jacobi's identity,
$$\aligned
q\psq 6\prod_{k=1}^{\infty}(1-q^{2k})^3&=q\Big(\sum_{m=0}^{\infty}
q^{3m(m+1)}\Big)\Big(\sum_{k=0}^{\infty}(-1)^k(2k+1)q^{k(k+1)}\Big)
\\&=\Big(\sum_{m=0}^{\infty}q^{3(2m+1)^2/4}\Big)
\Big(\sum_{k=0}^{\infty}(-1)^k(2k+1)q^{(2k+1)^2/4}\Big)
\\&=\sumn\Big(\sum\Sb k,m\in\{0,1,2,\ldots\}\\4n=(2k+1)^2+3(2m+1)^2
\endSb (-1)^k(2k+1)\Big)q^n\\&=\sumn\Big(\sum\Sb a,b\in\Bbb Z^+, 2\nmid
a\\4n=a^2+3b^2\endSb (-1)^{(a-1)/2}a\Big)q^n.
\endaligned$$
Hence
$$g_1(n)=\sum\Sb a,b\in\Bbb Z^+, 2\nmid
a\\4n=a^2+3b^2\endSb (-1)^{(a-1)/2}a.$$ It then follows that
$$\aligned &t(1,1,6,8;2n+1)=2t(1,3,4,4;n)
=2\Big(3^{\beta+1}-(-1)^{\beta}\Ls{-3}{n_1}\Big) \sum_{d\mid
n_1}\f{n_1}d\Ls 3d-4g_1(2n+3)
\\&=2\Big(3^{\beta+1}\Ls 3{n_1}-(-1)^{\beta}\Ls{-1}{n_1}\Big) \sum_{d\mid
n_1}d\Ls 3d-4\sum\Sb a,b\in\Bbb Z^+, 2\nmid a\\8n+12=a^2+3b^2\endSb
(-1)^{(a-1)/2}a.\endaligned$$ Observe that
$(-1)^{\beta}\ls{-1}{n_1}\e 3^{\beta}n_1=2n+3\e (-1)^{n-1}\mod 4$.
We then obtain part(ii). Now the proof is complete.

\section*{5. Special relations between $t(a,b,c,d;n)$ and \newline
$N(a,b,c,d;n)$}
\par In this section, we present some special relations between
 $t(a,b,c,d;n)$ and $N(a,b,c,d;n)$ and pose a lot of related
 conjectures.

 \pro{Theorem 5.1} Let
$n\in\Bbb Z^+$  Then
$$\align&t(1,4,7,8;n)=2N(1,4,7,8;2n+5)\qtq{for}n\e 3\mod 4,
\\&t(1,4,8,15;n)=2N(1,4,8,15;2n+7)\qtq{for}n\e 2\mod 4,
\\ &t(3,5,12,24;n)=2N(3,5,12,24;2n+11)\qtq{for}n\e 3\mod 4,
\\&t(3,5,20,40;n)=2N(3,5,20,40;2n+17)\qtq{for}n\e 3\mod 4.
\endalign$$
\endpro
Proof. From (3.15), (3.31) and (3.32) we see that
$$\align &t(1,4,7,8;4n+3)=t(1,2,7,4;2n+1)=t(1,1,2,7;n),
\\&t(1,4,8,15;4n+6)=t(1,1,2,15;n),
\\&t(3,5,12,24;4n+3)=t(3,3,5,6;n), \\&
t(3,5,20,40;4n+3)=t(3,5,5,10;n).\endalign$$
Now applying Theorem 2.2
yields the result.
\par\q
\pro{Theorem 5.2} For $n\in\Bbb Z^+$ we have
$$\align &t(2,3,3,4;n)=2N(2,3,3,4;2n+3)\qtq{for}n\e 2,3\mod 4,
\\&t(2,3,3,12;n)=2N(2,3,3,12;2n+5)\qtq{for}n\e 0,1\mod 4,
\\&t(2,3,3,24;n)=4N(2,3,3,24;2n+8)\qtq{for}n\e 2\mod 4,
\\&t(2,3,3,36;n)=2N(2,3,3,36;2n+11)\qtq{for}n\e 2,3\mod 4,
\\&t(1,1,6,12;n)=2N(1,1,6,12;2n+5)\qtq{for}n\e 0,3\mod 4,
\\&t(1,1,6,16;n)=\cases N(1,1,3,8;n+3)\qtq{if}n\e 2\mod 8,
\\4N(1,1,3,8;n+3)\qtq{if}n\e 4\mod 8.
\endcases
\endalign$$
\endpro
Proof. By Theorem 2.3 and (3.8),
$$N(2,3,3,4;8n+9)=t(1,2,3,3;n)=\f 12t(2,3,3,4;4n+3).$$
By (1.6) and (1.9),
$$\align &\sum_{n=0}^{\infty}N(2,3,3,4;n)q^n=\phq 2\phq 4\phq 3^2
\\&=(\phq 8+2q^2\psq{16})(\phq{16}+2q^4\psq{32})
\\&\q\times(\phq{24}^2+4q^{12}\psq{48}^2+4q^6\psq{24}^2+4q^3\phq{48}\psq{24}
+8q^{15}\psq{24}\psq{96}).
\endalign$$
Thus£¬
$$\align &\sum_{n=0}^{\infty}N(2,3,3,4;8n+7)q^{8n+7}
\\&=\phq 8\phq{16}\cdot 8q^{15}\psq{24}\psq{96}+2q^4\phq 8\psq{32}
\cdot 4q^3\phq{48}\psq{24}
\\&=8q^7\phq 8\psq{24}(\phq{48}\psq{32}+q^8\phq{16}\psq{96})
\\&=8q^7\phq 8\psq{24}\psq 8\psq{24}\endalign$$
and so
$$\sum_{n=0}^{\infty}N(2,3,3,4;8n+7)q^n=8\varphi(q)\psi(q)\psq
3^2.\tag 5.1$$
 On the other hand, using (1.5), (1.6) and (1.8) we
see that for $b\in\Bbb Z^+$,
$$\aligned &\sum_{n=0}^{\infty}t(2,3,3,4b;n)q^n
\\&=16\psq {2}\psq {3}^2\psq {4b}=16\phq {3}\psq {2}\psq {6}\psq {4b}
\\&=16(\phq {12}+2q^{3}\psq{24})(\phq{12}\psq {8}+q^{2}
\phq {4}\psq{24})\psq {4b}.\endaligned\tag 5.2$$ Thus
$$ \sum_{n=0}^{\infty}t(2,3,3,4b;4n+2)q^{4n+2}
=16q^{2}\phq {4}\psq {4b}\phq{12}\psq{24}$$ and so
$$\sum_{n=0}^{\infty}t(2,3,3,4b;4n+2)q^n=16\varphi(q)\psq
{b}\phq{3}\psq{6}=16\varphi(q)\psq 3^2\psq b.$$ This together with
(5.1) yields $t(2,3,3,4;4n+2)=2N(2,3,3,4;8n+7)$. Therefore
$t(2,3,3,4;n)=2N(2,3,3,4;2n+3)$ for $n\e 2,3\mod 4$. The remaining
results can be proved similarly.
\par\q
\pro{Lemma 5.1} Suppose $a,b\in\Bbb Z^+$ and $ab\e -1\mod 4$. Then
$$\align &\sumn N(a,a,a,2b;8n+5a)q^n=24\phq b\psq a\psq{2a}^2
\\&\sumn N(a,a,a,2b;8n+a+2b)q^n=12\phq a^2\psq a\psq{2b}.\endalign$$
\endpro
Proof. It is evident that
$$\align&\sumn N(a,a,a,2b;n)q^n
\\&=\phq a^3\phq{2b}=(\phq{16a}+2q^{4a}\psq{32a}+2q^a\psq{8a})^3
(\phq{8b}+2q^{2b}\psq{16b})
\\&=\big((\phq{16a}+2q^{4a}\psq{32a})^3+6q^a(\phq{16a}+2q^{4a}\psq{32a})^2\psq{8a}
+12q^{2a}(\phq{16a}\\&\q+2q^{4a}\psq{32a})\psq{8a}^2
+8q^{3a}\psq{8a}^3\big)(\phq{8b}+2q^{2b}\psq{16b}).
\endalign$$
Thus,
$$\align &\sumn N(a,a,a,2b;8n+5a)q^{8n+5a}
=6q^a\phq{16a}\cdot 4q^{4a}\psq{32a}\psq{8a}\phq{8b},
\\&\sumn N(a,a,a,2b;8n+a+2b)q^{8n+a+2b}
\\&\q=3(\phq{16a}^2+4q^{8a}\psq{32a}^2)\cdot 2q^a\psq{8a}\cdot
2q^{2b}\psq{16b}=12q^{a+2b}\phq{8a}^2\psq{8a}\psq{16b}
\endalign$$
and so the result follows.
\par\q

\pro{Theorem 5.3} For $n\in\Bbb Z^+$ with $n\e 3,5\mod 8$,
$$\align &t(1,1,2,12;n)=4N(1,1,4,6;n+2)=\f 83N(1,1,1,6;n+2)
\\&t(3,3,4,6;n)=\f 83N(2,3,3,3;n+2).\endalign$$
\endpro
Proof. By (3.1) and [S1, Theorem 2.11],
$$\align &t(1,1,2,12;8n+3)=2t(1,1,4,6;4n+1)=4t(1,2,3,4;2n)=4N(1,1,4,6;8n+5),
\\&t(1,1,2,12;8n+5)=2t(1,1,4,6;4n+2)=4N(1,1,4,6;8n+7).
\endalign$$
 From [S1, p.283] and Lemma 5.1 we know that
$$\align
&\sumn t(1,1,4,6;4n+1)q^n=32\phq 3\psi(q)\psq 2^2=\f 43\sumn
N(1,1,1,6;8n+5)q^n,
\\&\sumn t(1,1,4,6;4n+2)q^n=16\varphi(q)^2\psi(q)\psq 6
=\f 43\sumn N(1,1,1,6;8n+7)q^n.\endalign$$ Hence the formula for
$t(1,1,2,12;n)$ is true. From (5.2) we see that
$$\align &\sumn t(2,3,3,12;4n)q^{4n}=16\phq{12}^2\psq 8\psq{12},
\\&\sumn t(2,3,3,12;4n+5)q^{4n+5}=32q^5\psq{24}\phq
4\psq{24}\psq{12}.\endalign$$ Thus, appealing to Lemma 5.1 we get
$$\align &\sumn t(2,3,3,12;4n)q^n=16\phq{3}^2\psq 2\psq{3}
=\f 43\sumn N(2,3,3,3;8n+5)q^n,
\\&\sumn t(2,3,3,12;4n+5)q^n=32\varphi(q)\psq 3\psq 6^2
=\f 43\sumn N(2,3,3,3;8n+15)q^n.\endalign$$ This together with (3.1)
gives
$$\align &t(3,3,4,6;8n+3)=2t(2,3,3,12;4n)=\f 83\sumn N(2,3,3,3;8n+5),
\\&t(3,3,4,6;8n+5)=2t(2,3,3,12;4n+1)=4t(1,3,6,12;2n-1)\\&\qq\qq\qq\qq=
\f 83\sumn N(2,3,3,3;8n+15),
\endalign$$
which completes the proof.
\par\q
\pro{Corollary 5.1} For $n\in\Bbb Z^+$ with $n\e 5,7\mod 8$,
$$N(1,1,1,6;4n)=5N(1,1,1,6;n)\qtq{and}N(2,3,3,3;4n)=5N(2,3,3,3;n).$$
\endpro
Proof. From Theorem 5.3 and [S1, Theorem 2.1],
$$\align &\f 83N(1,1,1,6;n)=t(1,1,2,12;n-2)=\f 23(N(1,1,1,6;4n)-N(1,1,1,6;n)),
\\&\f 83N(2,3,3,3;n)=t(3,3,4,6;n-2)=\f 23(N(2,3,3,3;4n)-N(2,3,3,3;n)).\endalign$$
This yields the result.
\par\q
\par Using similar method one can prove the following results:
$$\align &t(3,3,4,18;n)=2N(3,3,4,18;2n+7)\qtq{for}n\e 0,1\mod 4,\tag 5.3
\\&t(1,3,8,12;n)=4N(1,3,8,12;n+3)\qtq{for}n\e 2,4\mod 8,\tag 5.4
\\&t(1,1,2,28;n)=4N(1,1,2,28;n+4)\qtq{for}n\e 1,3\mod 8,\tag 5.5
\\&t(1,3,4,24;n)=4N(1,3,4,24;n+4)\qtq{for}n\e 1,3\mod 8,\tag 5.6
\\&t(2,3,3,48;n)=N(2,3,3,48;2n+14)\qtq{for}n\e 0\mod 8,\tag 5.7
\\&t(1,1,8,14;n)=8N(1,1,8,14;n+3)\qtq{for}n\e 1\mod 8,\tag 5.8
\\&t(1,1,10,20;n)=4N(1,1,10,20;n+4)\qtq{for}n\e 1\mod 8,\tag 5.9
\\&t(1,1,14,16;n)=4N(1,1,14,16;n+4)\qtq{for}n\e 1\mod 8,\tag 5.10
\\&t(1,2,7,14;n)=8N(1,2,7,14;n+3)\qtq{for}n\e 1\mod 8,\tag 5.11
\\&t(1,1,8,30;n)=4N(1,1,8,30;2n+10)\qtq{for}n\e 1\mod 8,\tag 5.12
\\&t(1,3,4,16;n)=\f 43N(1,3,4,16;2n+6)\qtq{for}n\e 1\mod 8,\tag 5.13
\\&t(3,3,10,48;n)=4N(3,3,10,48;2n+16)\qtq{for}n\e 1\mod 8,\tag 5.14
\\&t(1,1,8,14;n)=4N(1,1,8,14;2n+6)\qtq{for}n\e 3\mod 8,\tag 5.15
\\&t(2,15,15,24;n)=4N(2,15,15,24;2n+14)\qtq{for}n\e 3\mod 8,\tag 5.16
\\&t(5,5,6,8;n)=4N(5,5,6,8;2n+6)\qtq{for}n\e 3\mod 8,\tag 5.17
\\&t(1,3,12,48;n)=\f 43N(1,3,12,48;2n+16)\qtq{for}n\e 4\mod 8,\tag 5.18
\\&t(2,4,5,5;n)=4N(2,4,5,5;n+2)\qtq{for}n\e 5\mod 8,\tag 5.19
\\&t(4,7,7,14;n)=4N(4,7,7,14;n+4)\qtq{for}n\e 5\mod 8,\tag 5.20
\\&t(1,1,16,30;n)=4N(1,1,16,30;2n+12)\qtq{for}n\e 5\mod 8,\tag 5.21
\\&t(1,1,30,40;n)=4N(1,1,30,40;2n+18)\qtq{for}n\e 5\mod 8,\tag 5.22
\\&t(1,3,16,36;n)=\f 43N(1,3,16,36;2n+14)\qtq{for}n\e 5\mod 8,\tag 5.23
\\&t(2,3,3,32;n)=4N(2,3,3,32;2n+10)\qtq{for}n\e 5\mod 8,\tag 5.24
\\&t(2,7,7,24;n)=4N(2,7,7,24;2n+10)\qtq{for}n\e 5\mod 8,\tag 5.25
\\&t(3,3,10,24;n)=4N(3,3,10,24;2n+10)\qtq{for}n\e 5\mod 8,\tag 5.26
\\&t(1,7,16,16;n)=4N(1,7,16,16;n+5)\qtq{for}n\e 6\mod 8,\tag 5.27
\\&t(2,3,3,48;n)=4N(2,3,3,48;2n+14)\qtq{for}n\e 6\mod 8,\tag 5.28
\\&t(1,1,10,20;n)=4N(1,1,10,20;n+4)\qtq{for}n\e 7\mod 8,\tag 5.29
\\&t(2,4,5,5;n)=4N(2,4,5,5;n+2)\qtq{for}n\e 7\mod 8,\tag 5.30
\\&t(4,7,7,14;n)=4N(4,7,7,14;n+4)\qtq{for}n\e 7\mod 8,\tag 5.31
\\&t(1,1,14,16;n)=4N(1,1,14,16;2n+8)\qtq{for}n\e 7\mod 8,\tag 5.32
\\&t(5,5,6,40;n)=4N(5,5,6,40;2n+14)\qtq{for}n\e 7\mod 8.\tag 5.33
\endalign$$
\par By doing calculations on Maple, we pose the following conjectures.

 \pro{Conjecture 5.1} Suppose $n\in\Bbb Z^+$. Then
$$\align&t(1,3,5,15;n)=8N(1,3,5,15;n+3)\q\t{for $n\e 3\mod 4$,}
\\&t(1,3,7,21;n)=8N(1,3,7,21;n+4)\q\t{for $n\e 2\mod 4$,}
\\&t(1,3,9,27;n)=\f{16}{13}N(1,3,9,27;2n+10)\q\t{for $n\e 7\mod 8$}.
\endalign$$\endpro
\par We remark that Dongxi Ye informed the author he proved the
formula for $t(1,3,5,15;n)$ and $t(1,3,7,21;n)$ by using theta
function identities and weight 2 modular forms.

\pro{Conjecture 5.2} Let $n\in\Bbb Z^+$. Then
$$\align&t(1,2,3,10;n)\\&\q=\cases\f 43N(1,2,3,10;2n+4)=\f
49N(1,2,3,10;8n+16)&\t{if $4\mid n-1$},
\\\f 83N(1,2,3,10;2n+4)&\t{if $8\mid n$,}
\\\f{16}9N(1,2,3,10;2n+4)&\t{if $n\e 10\mod{16}$,}
\\4N(1,2,3,10;2n+4)&\t{if $n\e 11,15\mod{20}$.}
\endcases\endalign$$\endpro

\pro{Conjecture 5.3} Let $n\in\Bbb Z^+$. Then
$$\align
&t(1,2,3,18;n)=\cases\f 43N(1,2,3,18;2n+6)=\f
49N(1,2,3,18;8n+24)&\t{if $4\mid n$},
\\\f 83N(1,2,3,18;2n+6)&\t{if $8\mid n-3$,}
\\ 4N(1,2,3,18;2n+6)&\t{if $12\mid n-6$,}
\\\f 85N(1,2,3,18;2n+6)&\t{if $24\mid n-15$,}
\endcases\endalign$$\endpro

\pro{Conjecture 5.4} Let $n\in\Bbb Z^+$. Then
$$\align
&t(1,3,6,30;n)\\&\q =\cases \f 43N(1,3,6,30;2n+10)=\f
49N(1,3,6,30;8n+40)&\t{if $4\mid n$,}
\\\f 83N(1,3,6,30;2n+10)=\f 8{15}N(1,3,6,30;8n+40)&\t{if $8\mid n-1$,}
\\\f{16}9N(1,3,6,30;2n+10)=\f{16}{33}N(1,3,6,30;8n+40)&\t{if $16\mid n+1$.}
\endcases\endalign$$\endpro

\pro{Conjecture 5.5} Let $n\in\Bbb Z^+$. Then
$$\align&t(1,3,18,18;n)\\&\q=\cases \f 43N(1,3,18,18;2n+10)=\f
49N(1,3,18,18;8n+40)&\t{if $4\mid n-2$,}
\\\f 83N(1,3,18,18;2n+10)=\f 8{15}N(1,3,18,18;8n+40)&\t{if $8\mid n-1$,}
\\\f {16}9N(1,3,18,18;2n+10)=\f {16}{33}N(1,3,18,18;8n+40)&\t{if $16\mid n-7$,}
\\ 4N(1,3,18,18;2n+10)=\f {4}{7}N(1,3,18,18;8n+40)&\t{if $12\mid n-4$,}
\\\f 85N(1,3,18,18;2n+10)=\f {8}{17}N(1,3,18,18;8n+40)&\t{if $24\mid n-13$.}
\endcases\endalign$$
\endpro

\pro{Conjecture 5.6} Let $n\in\Bbb Z^+$. Then
$$t(2,3,9,18;n)=\cases
\f 43N(2,3,9,18;2n+8)=\f 49N(2,3,9,18;8n+32)&\t{if $4\mid n-3$,}
\\\f 83N(2,3,9,18;2n+8)&\t{if $8\mid n-2$,}
\\\f {16}9N(2,3,9,18;2n+8)&\t{if $16\mid n-8$,}
\\\f {32}{15}N(2,3,9,18;2n+8)&\t{if $32\mid n-20$,}
\\\f 85N(2,3,9,18;2n+8)&\t{if $24\mid n-14$,}
\\4N(2,3,9,18;2n+8)&\t{if $12\mid n-5$.}
\endcases$$
\endpro

\pro{Conjecture 5.7} Let $n\in\Bbb Z^+$. Then
$$\align &t(2,5,10,15;n)\\&=\cases
\f 43N(2,5,10,15;2n+8)=\f 49N(2,5,10,15;8n+32)&\t{if $4\mid n-3$,}
\\\f 83N(2,5,10,15;2n+8)&\t{if $8\mid n-6$,}
\\\f{16}9N(2,5,10,15;2n+8)&\t{if $16\mid n-8$,}
\\4N(2,5,10,15;2n+8)&\t{if $n\e 61,81\mod{100}$}
\endcases\endalign$$\endpro

\pro{Conjecture 5.8} Let $n\in\Bbb Z^+$. Then
$$\align&t(5,6,15,30;n)=\cases \f 43N(5,6,15,30;2n+14)&\t{if $4\mid n-2$,}
\\\f 83N(5,6,15,30;2n+14)&\t{if $8\mid n-7$},
\\\f {16}9N(5,6,15,30;2n+14)&\t{if $16\mid n-13$}.
\endcases
\endalign$$\endpro

\pro{Conjecture 5.9} Let $n\in\Bbb Z^+$. Then
$$t(1,6,9,12;n)=2N(1,6,9,12;2n+7)\q\t{for $n\e 2,3\mod 4$}.$$
\endpro

\pro{Conjecture 5.10} Let $n\in\Bbb Z^+$ with $n\e 0\mod 8$. Then
$$\align &t(1,3,24,28;n)=4N(1,3,24,28;n+7),
\\&t(1,6,10,15;n)=8N(1,6,10,15;n+4),
\\&t(1,2,21,24;n)=4N(1,2,21,24;2n+12),
\\&t(1,6,8,33;n)=4N(1,6,8,33;2n+12).
\endalign$$\endpro

\pro{Conjecture 5.11} Let $n\in\Bbb Z^+$ with $n\e 1\mod 8$. Then
$$\align
&t(1,2,5,32;n)=4N(1,2,5,32;2n+10),
\\&t(1,2,13,24;n)=4N(1,2,13,24;2n+10),
\\&t(2,5,9,24;n)=4N(2,5,9,24;2n+10).
\endalign$$
\endpro
\pro{Conjecture 5.12} Let $n\in\Bbb Z^+$ with $n\e 2\mod 8$. Then
$$\align &t(1,2,5,8;n)=\f 43N(1,2,5,8;2n+4),
\\&t(1,2,5,24;n)=4N(1,2,5,24;2n+8).
\endalign$$
\endpro
\pro{Conjecture 5.13} Let $n\in\Bbb Z^+$ with $n\e 3\mod 8$. Then
$$\align&t(1,6,8,9;n)=4N(1,6,8,9;2n+6),
\\&t(1,18,24,45;n)=4N(1,18,24,45;2n+22),
\\&t(2,3,24,27;n)=4N(2,3,24,27;2n+14),
\\&t(2,5,10,15;n)=\f 43N(2,5,10,15;2n+8),
\\&t(3,6,8,39;n)=4N(3,6,8,39;2n+14).
\endalign$$\endpro

\pro{Conjecture 5.14} Let $n\in\Bbb Z^+$ with $n\e 4\mod 8$. Then
$$\align&t(1,8,16,31;n)=4N(1,8,16,31;n+7),
\\&t(1,9,30,40;n)=4N(1,9,30,40;2n+20),
\\&t(2,2,3,9;n)=\f 83N(2,2,3,9;2n+4).\endalign$$\endpro

\pro{Conjecture 5.15} Let $n\in\Bbb Z^+$ with $n\e 5\mod 8$. Then
$$\align &t(1,8,16,23;n)=4N(1,8,16,23;n+6),
\\&t(1,5,10,40;n)=\f 43N(1,5,10,40;2n+14),
\\&t(2,3,11,24;n)=4N(2,3,11,24;2n+10).
\endalign$$\endpro

\pro{Conjecture 5.16} Let $n\in\Bbb Z^+$ with $n\e 6\mod 8$. Then
$$\align &t(1,18,21,24;n)=4N(1,18,21,24;2n+16),
\\&t(3,6,8,15;n)=4N(3,6,8,15;2n+8),
\\&t(5,6,8,45;n)=4N(5,6,8,45;2n+16).\endalign$$
\endpro

\pro{Conjecture 5.17} Let $n\in\Bbb Z^+$ with $n\e 7\mod 8$. Then
$$\align &t(2,3,5,30;n)=8N(2,3,5,30;n+5),
\\&t(1,6,9,24;n)=\f 43N(1,6,9,24;2n+10),
\\&t(2,5,24,25;n)=4N(2,5,24,25;2n+14),
\\&t(2,9,21,24;n)=4N(2,9,21,24;2n+14),
\\&t(3,6,7,8;n)=4N(3,6,7,8;2n+6),
\\&t(3,15,30,40;n)=4N(3,15,30,40;2n+22).
\endalign$$
\endpro
\pro{Conjecture 5.18} Let $n\in\Bbb Z^+$. Then
$$\align &t(1,1,25,25;n)=2N(1,1,25,25;2n+13)\qtq{for}n\e 0,2\mod 5,
\\&t(1,1,5,21;n)=2N(1,1,5,21;2n+7)\qtq{for}n\e 0,3\mod 5,
\\&t(1,6,6,15;n)=2N(1,6,6,15;2n+7)\qtq{for}n\e 0,3\mod 5,\\&t(2,2,3,5;n)=2N(2,2,3,5;2n+3)\qtq{for}n\e 3,4\mod 5,
\\&t(3,7,7,35;n)=2N(3,7,7,35;2n+13)\qtq{for}n\e 3,4\mod 5,
\\&t(3,3,7,15;n)=2N(3,3,7,15;2n+7)\qtq{for}n\e 1,2\mod 5.
\\ &t(1,1,3,5;n)=\f 25N(1,1,3,5;8n+10)\qtq{for}n\e 2,3\mod 5,
\\&t(1,3,3,15;n)=\f 25N(1,3,3,15;8n+22)\qtq{for}n\e 0,2\mod 5.
\endalign$$
\endpro

\pro{Conjecture 5.19} Let $n\in\Bbb Z^+$ with $n\e 0,1,5\mod 7$.
Then
$$t(1,2,2,7;n)=2N(1,2,2,7;2n+3).$$
\endpro
\pro{Conjecture 5.20} Let $n\in\Bbb Z^+$ with $n\e 0,2,3\mod 7$.
Then
$$t(1,1,1,7;n)=\f 25N(1,1,1,7;8n+10).$$
\endpro

\pro{Conjecture 5.21} Let $n\in\Bbb Z^+$. Then
$$\align &t(1,1,1,33;n)=2N(1,1,1,33;2n+9)\qtq{for}n\e 2,4,5,6,10\mod{11},
\\&t(1,1,9,33;n)=2N(1,1,9,33;2n+11)\qtq{for}n\e 1,3,4,5,9\mod{11},
\\&t(1,9,9,33;n)=2N(1,9,9,33;2n+13)\qtq{for}n\e 0,2,3,4,8\mod{11}.
\endalign$$
\endpro
\par{\bf Remark 5.1} From [S3, Theorem 4.2] we know that if $a,b,c,d,n\in\Bbb Z^+$ and
$a\e b\e c\e d\e \pm 1\pmod 4$, then
$$t(a,b,c,d;n)=N\big(a,b,c,d;8n+a+b+c+d\big)-
N(a,b,c,d;2n+(a+b+c+d)/4).$$


\begin{thebibliography}{BCH}
\bibitem [ACH]{} C. Adiga, S. Cooper and J. H. Han, {\it
A general relation between sums of squares and sums of triangular
numbers}, Int. J. Number Theory {\bf 1}(2005), 175-182.

\bibitem [BCH]{} N. D. Baruah, S. Cooper and M. Hirschhorn, {\it
Sums of squares and sums of triangular numbers induced by partitions
of $8$}, Int. J. Number Theory {\bf 4}(2008), 525-538.

\bibitem [Be] {} B.C. Berndt, {\it Ramanujan's Notebooks}, Part III,
Springer, New York, 1991.



\bibitem [C]{} S. Cooper, {\it On the number of
representations of integers by certain quadratic forms, II}, J.
Combin. Number Theory {\bf 1}(2009), 153-182.

\bibitem [D]{} L.E. Dickson, {\it History of the Theory of Numbers},
 Vol. II, Carnegie Institute of Washington,
Washington D.C., 1923. Reprinted by AMS Chelsea, 1999.

\bibitem [S1] {} Z.H. Sun, {\it Some relations between
$t(a,b,c,d;n)$ and $N(a,b,c,d;n)$}, Acta Arith. {\bf 175}(2016),
269-289.

\bibitem [S2] {} Z.H. Sun, {\it Ramanujan's theta functions and sums
of triangular numbers}, Int. J. Number Theory {\bf 15}(2019),
969-989.

\bibitem [S3] {} Z.H. Sun, {\it The number of representations of
$n$ as a linear combination of triangular numbers}, Int. J. Number
Theory {\bf 15}(2019), 1191-1218.

\bibitem [XZ] {} E.X.W. Xia and Y. Zhang, {\it
Proofs of some conjectures of Sun on the relations between sums of
squares and sums of triangular numbers}, Int. J. Number Theory {\bf
15}(2019), 189-212.

\bibitem [Y]{} X.M. Yao, {\it The relations between $N(a,b,c,d;n)$
and $t(a,b,c,d;n)$ and $(p,k)$-parametrization of theta functions},
J. Math. Anal. Appl. {\bf 453}(2017), 125-143.

\bibitem [WS1] {} M. Wang and Z.H. Sun, On the number of
representations of $n$ as a linear combination of four triangular
numbers, Int. J. Number Theory {\bf 12}(2016), 1641-1662.

\bibitem [WS2] {} M. Wang and Z.H. Sun, {\it On  the number of
representations of $n$ as a linear combination of four triangular
numbers II}, Int. J. Number Theory {\bf 13}(2017), 593-617.

\bibitem [W] {} K.S. Williams, {\it Number Theory in the Spirit of
Liouville}, Cambridge Univ. Press, New York, 2011.

\end{thebibliography}
\end{document}